\documentclass[11pt]{amsart}
\usepackage[english]{babel}

\usepackage{latexsym}
\usepackage{amssymb}
\usepackage{amsmath}
\usepackage{amsfonts}
\usepackage[mathscr]{eucal}
\usepackage{color}
\usepackage{setspace}
\usepackage{algorithm}
\usepackage{algpseudocode}

\usepackage{amscd}

\usepackage{graphicx}
\usepackage{epsfig}
\usepackage{relsize}

\usepackage{fullpage}

\usepackage[table,xcdraw]{xcolor}
\usepackage{longtable}
\usepackage{rotating}
\usepackage{array}
\usepackage{multicol}
\usepackage{multirow}
\usepackage{makecell}
\usepackage{mathtools}
\usepackage{booktabs}

\usepackage{amsthm}
\usepackage[all,cmtip]{xy}
\usepackage[pdfstartview=FitH]{hyperref}
\usepackage{comment}
\usepackage{psfrag}
\usepackage{verbatim}
\usepackage{xurl}
\hypersetup{
colorlinks=true,
linkcolor=blue
}
\usepackage{pifont}
\usepackage{caption}
\usepackage{csquotes}
\usepackage{lscape}
\usepackage{pdflscape}
\usepackage{everypage}
\usepackage{afterpage}
\usepackage{arydshln}

\usepackage{soul}
 \usepackage[normalem]{ulem}

\usepackage{tikz}
\usepackage{pgfplots}
\pgfplotsset{compat=1.15}
\usepackage{mathrsfs}
\usetikzlibrary{arrows}
\usetikzlibrary{patterns}
\usetikzlibrary{calc}
\usetikzlibrary{decorations.markings,positioning}

\usepackage[style=numeric,backend=biber,maxnames=8]{biblatex}

\newtheorem*{thm*}{Theorem}
\newtheorem*{cor*}{Corollary}		

\theoremstyle{definition}
\newtheorem*{rema}{Remark} 
\newtheorem*{definition*}{Definition}    	

\newtheorem*{ack}{Acknowledgements}

\numberwithin{equation}{section}

\DeclareMathOperator{\dist}{dist}

\newcommand{\norm}[1]{\left\Vert#1\right\Vert}

\newcommand{\RR}{\mathbb{R}}

\def\boxit#1#2{%
    \smash{\color{red}\fboxrule=1pt\relax\fboxsep=2pt\relax%
    \llap{\rlap{\fbox{\phantom{\rule{#1}{#2}}}}~}}\ignorespaces
}

\def\boxittwo#1#2{%
    \smash{\color{blue}\fboxrule=1pt\relax\fboxsep=2pt\relax%
    \llap{\rlap{\fbox{\phantom{\rule{#1}{#2}}}}~}}\ignorespaces
}

\def\boxitthree#1#2{%
    \smash{\color{orange}\fboxrule=1pt\relax\fboxsep=2pt\relax%
    \llap{\rlap{\fbox{\phantom{\rule{#1}{#2}}}}~}}\ignorespaces
}

\begin{document}
























\title[TDA in ATM]{Topological Data Analysis in Air Traffic Management: the shape of big flight data sets \\ 
}

\author[M.~Cuerno]{Manuel Cuerno$^{\ast}$}

\author[L.~Guijarro]{Luis Guijarro$^{\ast\ast}$}

\author[R.~Arnaldo]{Rosa María  Arnaldo Valdés$^{\ast\ast\ast}$}

\author[F.~Gómez]{Víctor Fernando  Gómez Comendador$^{\ast\ast\ast\ast}$}

\thanks{$^{\ast}$Supported in part by the FPI Graduate Research Grant PRE2018-084109, and by research grants  
	 MTM2017-‐85934-‐C3-‐2-‐P and PID2021-124195NB-C32
from the Ministerio de Econom\'ia y Competitividad de Espa\~{na} (MINECO)}
\thanks{$^{\ast\ast}$Supported by research grants  MTM2017-‐85934-‐C3-‐2-‐P and PID2021-124195NB-C32 from the  Ministerio de Econom\'{\i}a y Competitividad de Espa\~{na} (MINECO), and by ICMAT Severo Ochoa project SEV-‐2015-‐0554 (MINECO)} 
\thanks{$^{\ast\ast\ast}$Universidad Politécinca of Madrid (UPM)}
\thanks{$^{\ast\ast\ast\ast}$Universidad Politécnica of Madrid (UPM)}


\address[M.~Cuerno]{Department of Mathematics, CUNEF, Spain}
\email{manuel.mellado@cunef.edu}

\address[L.~Guijarro]{Department of Mathematics, Universidad Aut\'onoma de Madrid and ICMAT CSIC-UAM-UC3M, Spain}
\email{luis.guijarro@uam.es} 

\address[R.~Arnaldo]{Department of aerospace systems, air transport and airports, Universidad Polit\'ecnica de Madrid, Spain}
\email{rosamaria.arnaldo@upm.es}

\address[F.~Gómez]{Department of aerospace systems, air transport and airports, Universidad Polit\'ecnica de Madrid, Spain}
\email{fernando.gcomendador@upm.es}


\date{\today}





\begin{abstract}

\textcolor{black}{Analyzing flight trajectory data sets poses challenges due to the intricate interconnections among various factors and the high dimensionality of the data. Topological Data Analysis (TDA) is a way of analyzing big data sets focusing on the topological features this data sets have as point clouds in some metric space. Techniques as the ones that TDA provides are suitable for dealing with high dimensionality and intricate interconnections. This paper introduces TDA and its tools and methods as a way to derive meaningful insights from ATM data. Our focus is on employing TDA to extract valuable information related to airports. Specifically, by utilizing persistence landscapes (a potent TDA tool) we generate footprints for each airport. These footprints, obtained by averaging over a specific time period, are based on the deviation of trajectories and delays. We apply this method to the set of Spanish' airports in the Summer Season of 2018. Remarkably, our results align with the established Spanish airport classification and raise intriguing questions for further exploration. This analysis serves as a proof of concept, showcasing the potential application of TDA in the ATM field. While previous works have outlined the general applicability of TDA in aviation, this paper marks the first comprehensive application of TDA to a substantial volume of ATM data. Finally, we present conclusions and guidelines to address future challenges in the ATM domain.}

\end{abstract}
\setcounter{tocdepth}{1}

\maketitle






\section{Introduction}

Airports and Air Traffic Management Systems form intricate sociotechnical structures, presenting challenges in analysis due to their high interdependence. Globally, these systems manage operations involving over 32,000 in-service aircraft operated by more than 1,300 airlines, facilitating over 4.1 billion passengers across 41.9 million flights and 45,000 routes at 3,700 airports \cite{uno,dos}.

The interconnectedness, interdependences and complexity of the system is reflected in the amount of data generated by its operation. Although flight trajectory data offer a big potential to grasp the features and behaviour of such complex system, analyzing flight trajectory data proves challenging due to its high dimensionality, continuous nature, and multi-level interactions influenced by factors such as air traffic control, weather conditions, aircraft type, and pilot behaviour \cite{tres}. The dynamic changes and multiple variables involved, like altitude, speed, heading, longitude, latitude, and time, further complicate manual processing \cite{cuatro,seis,cinco}. Its high dimensionality can make it difficult to visualize or analyse the data without reducing its dimensionality through techniques like principal component analysis or t-SNE.

To address the high dimensionality issue, dimensionality reduction techniques such as principal component analysis (PCA) and t-distributed stochastic neighbour embedding (t-SNE) are employed to identify common traffic patterns \cite{diez,nueve,doce,ocho,catorce,once,trece}. However, limitations exist, including empirical selection of techniques and potential loss of relevant information, potential lack of interpretability, and issues of generalizability to different contexts. The choice between PCA or t-SNE in the analysis algorithm is often made through empirical observations, and alternative machine learning techniques might yield superior results for specific data types or traffic patterns. While PCA and t-SNE can uncover data patterns, the interpretation of these patterns poses challenges, requiring additional domain expertise for meaningful insights. The studies conducted may lack generalizability to different airports or contexts, necessitating further research to assess the applicability of these techniques in diverse settings. Moreover, some studies encounter issues related to limited sample size, data quality, or pre-processing, and there's always a risk of losing relevant information during dimensional reduction.

Clustering techniques like k-means, hierarchical clustering, and DBSCAN are alternative methods to identify common traffic patterns \cite{ocho}. In \cite{dieciseis}, a trajectory clustering ensemble method is introduced, utilizing a similarity matrix and employing the Nanjing Lukou Airport terminal area as an illustrative example. In the study by Zeng et al \cite{diecisiete}, a DBSCAN clustering analysis method is employed for the detection of outliers in aircraft airborne and controllers' data, contributing to the assessment and monitoring of flight status. Furthermore, a paper at the SESAR Innovation Days 2017 \cite{dieciocho} proposed a method for analyzing the traffic patterns of civil aviation flights using principal component analysis (PCA) and DBSCAN clustering. The method was tested on a dataset of flight trajectory data from the Beijing Capital International Airport, and the results showed that the method was able to identify common traffic patterns and distinguish them from those abnormal.

While effective, they have limitations such as sensitivity to initialization, difficulty in determining the number of clusters and handling noise, scalability and interpretability issues and outlier sensitivity. Additionally these methods are limited to Euclidean distances (k-means and hierarchical clustering), to finding geometric structures and to unsupervised learning.

Beyond the previously mentioned methods, machine learning algorithms, including decision trees, random forests, and neural networks, are also utilized for pattern identification in flight trajectory data. However, limitations include sensitivity to data availability, lack of standardization, and challenges in result interpretation \cite{zhourobustness}.
As can be seen, overall, while the existing research works have made significant contributions to the field of analysing high-dimensional flight trajectory data, there are still limitations that need to be addressed in future research to improve the accuracy, reliability, and scalability of these methods. 

To overcome existing limitations, this paper proposes the use of Topological Data Analysis (TDA), a mathematical framework using tools from algebraic topology. TDA offers advantages in identifying topological features, simplifying complex data, handling noise, comparing datasets, and building predictive models \cite{chazal}.

A key strength of TDA lies in its ability to identify underlying structures and patterns in data with enhanced resilience to noise and outliers. For instance, Casacuberta et al. \cite{casacuberta} demonstrated this by estimating the intrinsic dimension of a point cloud in a high-dimensional space, revealing that the real dimension of the data was less than that of the ambient space. TDA also proves beneficial in mitigating the challenge of limited data availability, enabling the integration of diverse data sources and extracting meaningful insights from incomplete or noisy datasets. Furthermore, TDA offers a solution to scalability issues by facilitating the analysis of large and complex datasets through parallel computing techniques. Additionally, TDA addresses concerns related to the lack of standardization by providing a flexible and adaptable framework applicable to a broad spectrum of data formats and types, and  it works no matter the initial domain is and can give conclusions only based on the data input no matter where this data comes from. This allows the researcher extract some initial conclusions excluding the domain of the data.

TDA applications include identifying topological features, simplifying complex data, handling noise, comparing datasets, and building predictive models \cite{casacuberta,maxdetection}. In addition, other authors have also applied different topological tools for air traffic management problems. Examples include detecting anomalous trajectories and assessing delay network robustness in adverse weather conditions \cite{zhourobustness}, predictive modelling \cite{chourbanairspace}  or more topological graph-based studies related to finding patterns in different air traffic situations \cite{topological_characteristics}. 

Overall, TDA offers a powerful tool, either alone or combined with other existent methods, for analysing complex and high-dimensional aviation data sets. By identifying topological features and patterns, TDA can reveal hidden relationships and help airlines and airports making better decisions about flight scheduling, maintenance, and safety as well as improving airport operations and the passenger experience. 

In this work, we present a proof of concept of TDA applied to the Air Traffic and Transportation Management field. In particular, we construct a point cloud based on spatiotemporal and delay flight-data of a certain airport in order to develop, using TDA techniques, a footprint of it. After this process, we compare the actual classification of the Spanish Airport net (only based on number of passengers per year) with the one we produce with the footprints  and compute a \textit{block permutation test} in order to calculate the statistical significance of our experiment. Additionally, we compare our results with ones produced by a non-TDA method based on centrality measures.  Interesting similarities and discrepancies appeared on this experiment, showing how TDA can enrich the study of the ATM in future projects. 

To the best of our knowledge, there have been generic works outlining the possible application of TDA in aviation \cite{duponchel2018remote,maxtda} but none of them have applied any TDA technique apart from computing some simplicial complexes based on airports \cite{maxtda} or some sensor analysis partially related with aircraft search \cite{duponchel2018remote}. No rigorous work has already been conducted to apply TDA to a large volume of aircraft trajectory data for the purpose of anticipating and identifying deviations and anomalies in aircraft space/time trajectories. 

Additionally, there is a lack of research aimed at inferring patterns of behaviour at different airports or classifying and characterizing airports based on the distribution of their daily flights through trajectory deviation and delay.  Despite this gap in the literature, some works based on centrality measures and graph theory have been developed to globally classify networks of airports \cite{huynh2024understanding,nikolaou2020identification,song2017analysis,sun2023centrality}.

The application of TDA techniques, and other types of mathematical analysis algorithms applied to data series, requires close collaboration between research groups focused on mathematical modelling and analysis, and research groups focused on different fields of application. This means that certain lines identified as potentially of interest may lack continuity. In this case, other references have been identified that show how some groups, which have initiated work in these areas, maintain interest in this application and continue to develop various projects specifically \cite{frewin2023aggregate, holdren2023adaptable,li2021identification,li2021dynamics, li2021graph}. Moreover, TDA and other topological tools are also being applied in other areas of transport, not only in aviation, that surely could improve its implementation in the ATM field \cite{beutenmuller2023topological,duponchel2018remote,macarasig2024topological,zhang2024topological}.

This paper proposes a method for assessing the structural characteristics of air traffic situation based on TDA, providing new clues to give a more precise description of air traffic complexity, to assess the deviation from expected aircraft trajectories, to study the generation of delays, to identify structural patterns, and to detect and analyze anomalies in airport operations. In particular this method provides a natural classification of airports attending to operational characteristics. The method adjusts and captures, from a mathematical point of view, the grouping criteria that the airport provider uses operationally. It mimics those criteria. It has the potential to capture more complex operational criteria and relationships.

Furthermore, this initial work wants to encourage the future use of TDA in the aerospace field as a powerful tool to deal with the existent data constrains and enrich its combination with already known methods obtaining the optimal use of all of them.

\begin{ack}

The authors want to thank Telmo Pérez Izquierdo and Juan Antonio Pichardo for their comments and ideas to the manuscript. They have undoubtedly improved the quality of the material we present.

\end{ack}

\section{Material and methods}{\label{methods}}

Topological Data Analysis \cite{carlsson, chazal, otter2017roadmap} employs a set of structures, named simplicial complexes, for the purpose of data analysis. Roughly speaking, if the data is represented by a point cloud, we can endow it with a simplicial structure and study how it evolves when we vary certain scaling parameters. This is the setting of the one-parameter TDA, which is the one most commonly used, and the one we are going to implement in this work.

While multi-parameter TDA is an attractive field with promising results \cite{botnanmulti, harringtonmulti}, it is computationally challenging and thus most approaches are theoretical. However, the study of multi-parameter TDA could be an interesting setting to further explore problems in this field allowing us to include more factors in the analysis, as we have explained in the Introduction of this paper.

The primary objective of this Section is to offer introductory concepts of one-parameter TDA and the tools and methods we use in this paper. For the sake of clarity and concision, we will not present a more formal description of simplicial complexes, homology and persistent homology.  For the interested reader we recommend \cite{carlsson,chazal} and, a mathematical introduction to homology can be found in \cite{hatcher,weibel}.

\subsection{Persistent homology}

We will consider a \textit{point cloud} $\mathcal{V}$ as a set of points in $\RR^d$ for some $d\geq 2$. Roughly speaking, TDA would extract the most relevant ``topological features'', called $n$-cycles, of $\mathcal{V}$. The $n$-cycles represent $n$-dimensional holes. For example, $\mathbb{S}^1\subset\RR^2$ has one $1$-cycle. In order to obtain the information of the $n$-cycles in our point cloud, TDA usually computes the \textit{persistent homology} of $\mathcal{V}$. It will measure how predominant those $n$-cycles are.

To extract the persistent homology of the point cloud $\mathcal{V}$, we first need to construct a \textit{filtration} based on $\mathcal{V}$. A commonly used method for this is the \textit{Vietoris-Rips filtration}, denoted as $VR(\mathcal{V})$. This filtration is a nested sequence of simplicial complexes ($VR_{r_0}(\mathcal{V}) \subseteq VR_{r_1}(\mathcal{V})$ for $r_0 \leq r_1$, where $A \subseteq B$ means $A$ is a subset of $B$), and each complex $VR_r(\mathcal{V})$ in the filtration is parameterized by a radius $r > 0$.

The vertices (0-simplices) of these complexes correspond to the points in $\mathcal{V}$. A $k$-simplex, which can be understood as a $k$-dimensional tetrahedron, is included in $VR_r(\mathcal{V})$ if every pair of points in any $(k+1)$-tuple ${v_0, \dots, v_k} \subseteq \mathcal{V}$ satisfies $\dist_{\RR^d}(v_i, v_j) \leq r$ for all $i, j$. In other words, the collection of points ${v_0, \dots, v_k} \subseteq \mathcal{V}$ will form a $k$-dimensional tetrahedron. So, a simplicial complex consists of simplices of various dimensions.

\begin{rema}
Typically, the point clouds we deal with are finite, which means there exists some $r_{\max} \in \mathbb{R}^+$ such that $VR_r(\mathcal{V}) = VR_{r_{\max}}(\mathcal{V})$ for all $r \geq r_{\max}$.

Due to the construction of the filtration, all values of $r$ are in $[0,r_{\max}]$ are considered. Specifically, $VR_0(\mathcal{V}) = \mathcal{V}$, and for $r > r_{\max}$, we have $VR_r(\mathcal{V}) = VR_{r_{\max}}(\mathcal{V})$. Hence, no arbitrary choice of $r$ is made at any point in the process.

Thus, we obtain the following filtration \[
\mathcal{V}=VR_0(\mathcal{V})\subseteq\cdots\subseteq VR_r(\mathcal{V})\subseteq\cdots\subseteq VR_{r_0'}(\mathcal{V})=VR_{r_{\max}}(\mathcal{V}).
\]
\end{rema}

These inclusions induce maps $H_n(VR_{r_1}(\mathcal{V}))\to H_n(VR_{r_1}(\mathcal{V}))$ for $r_0\leq r_1$, where $H_n(VR_{r}(\mathcal{V}))$ is the $n$-homology group of $VR_r(\mathcal{V})$, i. e., the group spanned by the generators of each $n$--dimensional hole of $VR_r(\mathcal{V})$.  The \textit{$n$-persistent homology} of $\mathcal{V}$ with the Vietoris--Rips filtration measures how much the $n$-topological features live during the filtration. Every of those features has two coordinates $(b,d)$ that give the following information: $b$ (\textit{birth parameter}) is the smallest $r>0$ such that $H_n(VR_b(\mathcal{V}))$ has that feature and $d$ (\textit{death parameter}) is the biggest $r>b>0$ such that $H_n(VR_d(\mathcal{V}))$ also has it. If $d$ does not exists, we will assign $d=\infty$. The \textit{persistence} or \textit{lifetime} of an $n-$persistent cycle is defined as the difference $d-b$. For graphical intuition about these constructions, we recommend the beautiful figures made by R. Ghrist in the following paper \cite{ghrist}.

\subsection{Persistence diagrams and barcodes}

Once we know which information we are extracting for our point cloud, we need to represent those results. For that purpose, TDA mainly uses two tools: \textit{persistence diagrams} and \textit{barcodes}. These are powerful visual tools that keep all the information related with the persistent homology of our filtration.

There exist multiple equivalent definitions of  persistence diagrams. In order to mantain it as simple as possible, we will present them as a multiset of points in the upper semispace of the first quadrant of $\RR^2$,  i. e., $\RR^2_{\geq} = \{(x,y)\in\RR^2\colon y\geq x\}$. A \textit{multiset} $V$ of $\RR^d$ is a set of points in $\RR^d$ whose elements have multiplicity, i.e., each point in $V$ can be repeated more than once.
For example,  $V_0=\{(0,0),(1,0)\}$ and $V_1=\{(1,0),(1,0), (-1,1),(2,2), (2,2), (2,2)\}$ are multisets of $\RR^2$. In $V_0$ each element has multiplicity $1$ and in $V_1$, $(1,0)$ has multiplicity $2$ and $(2,2)$ has multiplicity $3$.

With this notion we are ready to present an easy definition of the \textit{persistence diagram} as in the spirit of  \cite{guijarrodiagram}.
A \textit{persistence diagram} $PD_i(\mathcal{F})$ is a multiset of points $(b,d)\in\overline{\mathbb{R}}^2_{\geq0}$ where $\overline{\mathbb{R}}^2_{\geq0}=\{(x,y)\in\overline{\mathbb{R}}\times\overline{\mathbb{R}}\colon 0\leq x<y\}$ and $\overline{\mathbb{R}}=\mathbb{R}\cup\{-\infty,\infty\}$, and  whose points represent the $i$--persistent homology of the filtration $\mathcal{F}$. Every point $(b,d)\in PD_i(\mathcal{F})$ encodes the information of an $i$--dimensional hole of the filtration $\mathcal{F}$ that appears at time $b$ and dies at time $d$.

\begin{figure}[h]
    \centering
    \includegraphics[width=430pt]{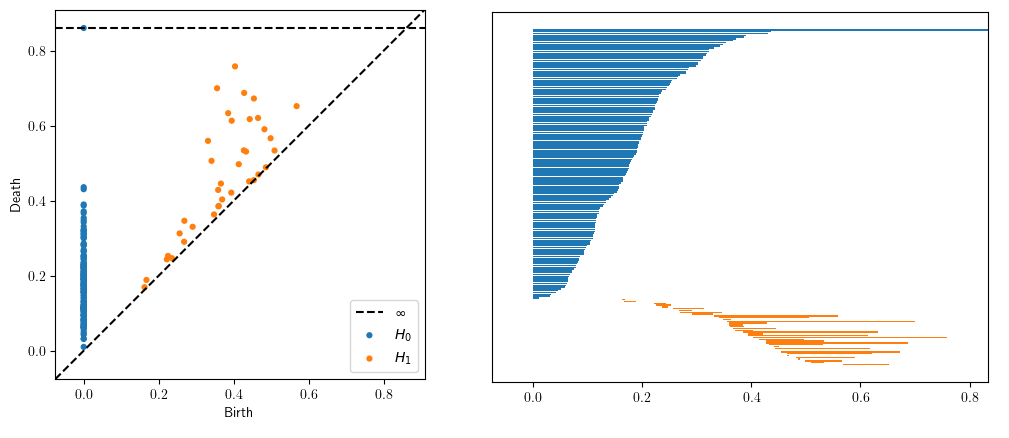}
    \caption{Persistence Diagram and Barcode of $H_0$ and $H_1$ for 150 random points on $\mathbb{S}^2\subset\mathbb{R}^3$. Regarding the left figure, the x-axis represents the birth and the y-axis the death of each topological feature. Regarding the right figure, the x-axis represents the birth and the death of each topological feature.}
    \label{PDbarcodeshpere}
\end{figure}

Equivalently, we can introduce the definition of a \textit{barcode} of certain persistence diagram $PD_i(\mathcal{F})$ as a set of intervals $[b, d]$, where each interval corresponds to a point $(b, d) \in PD_i(\mathcal{F})$.  Figure \ref{PDbarcodeshpere} illustrate that barcode allows to show the multiplicity of points whereas persistent diagrams not.

\begin{rema}
    In order to improve the readability of the following definitions and sections, we will use $PD_i(X)$ as the persistence diagram of the Vietoris--Rips filtration of the point cloud $X$.
\end{rema}


    Persistent diagrams and barcodes can be considered as visual tools to help in the interpretation of the persistent homology; we can combine various $i$--persistent homologies in the same picture, as we have done in Figure \ref{PDbarcodeshpere} with the $0$ and $1$--persistent homology of a point cloud of $150$ points extracted randomly from the unit sphere $\mathbb{S}^2\subset\mathbb{R}^3$.

    Finally, it is worth to mention that there exist different notion of distances between persistence diagrams such as the \textit{bottleneck} or the \textit{heat kernel} distance. All of those have their respective stability theorems which hold that small perturbations in the point clouds reflects on small perturbations on the persistence diagrams (please, consult \cite{chazalgromovhausdorff}).

\subsection{Persistence landscapes}

As discussed, persistence diagrams and barcodes serve as powerful tools for visualizing the outcomes of persistent homology. However, one limitation of these representations is that the information they offer is not amenable to straightforward algebraic manipulations. In other words, adding two barcodes or computing averages is not a well-defined process.

In order to deal with this problem, there exist several vectorizations of the persistent diagrams (see \cite{mariajose_survey}). In particular, in this paper we will use the \textit{persistence landscapes}. Bubenik and collaborators \cite{bubenikgradedandlandscapes,bubeniklandscapesstatistical,bubeniklandscapesproperties,bubeniklandscapestoolbox} defined them as functions that encapsulate the information of a given persistent diagram $PD_i(X)$. 



    Let $PD_i(X)=\{(a_j,b_j)\}_{j\in J}$ be a persistence diagram. We define the \textit{kth--persistence landscape} of $PD_i(X)$ via an auxiliary function: first, for $a<b$, let \[
        f_{(a,b)}(t)=\max(0,\min(t-a,b-t)).
        \]Then define \[
        \lambda_k(t)=\max{}_k\{f_{(a_j,b_j)}(t)\}_{j\in J},
        \]where $\max_k$ denotes the $k$-th largest element.
    
    This process is repeated until $\lambda_l\equiv 0$ for any $l\in\mathbb{N}$. The set of functions $\lambda = \{\lambda_1,\dots,\lambda_{l-1}\}$ is denoted as the \textit{persistence landscape} of $PD_i(X)$.


\begin{figure}[h]
    \centering
    \includegraphics[width=\columnwidth]{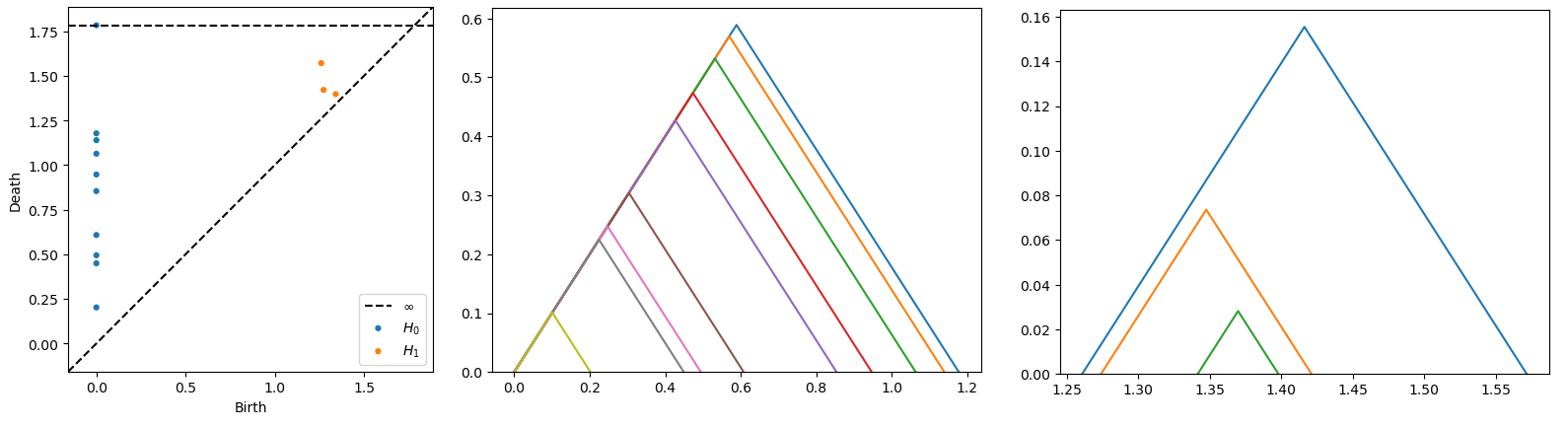}
    \caption{On the left: $PD_0$ and $PD_1$ of $10$ points on $\mathbb{S}^2\subset\mathbb{R}^3$. On the center: the landscapes of the $H_0$. On the right: the landscapes of the $H_1$.}
    \label{PDylandscapes}
\end{figure}


One of the advantages of persistence landscapes is the possibility of computing the \textit{average persistent landscape} in the presence of a family of persistent landscapes.
\color{black} Let $\lambda^{(1)},\lambda^{(2)},\dots,\lambda^{(n)}$ be a sequence of persistence landscapes with  their respective kth-persistence landscapes $\lambda^{(i)}=\{\lambda^{(i)}_1,\dots,\lambda^{(i)}_l\}$  for $i\in\{1,\dots,n\}$. Then we define the \textit{average of persistence landscapes} as \[
    \overline{\lambda}(k,t)=\frac1n\sum_{i=1}^n\lambda^{(i)}(k,t).
    \]In other words, each  kth-persistence landscape  of the new average persistence landscape is \[
    \overline{\lambda}_k(t)=\frac1n\sum_{i=1}^n\lambda_k^{(i)}(t).
    \]



\begin{figure}[h]
    \centering
    \includegraphics[width=\columnwidth]{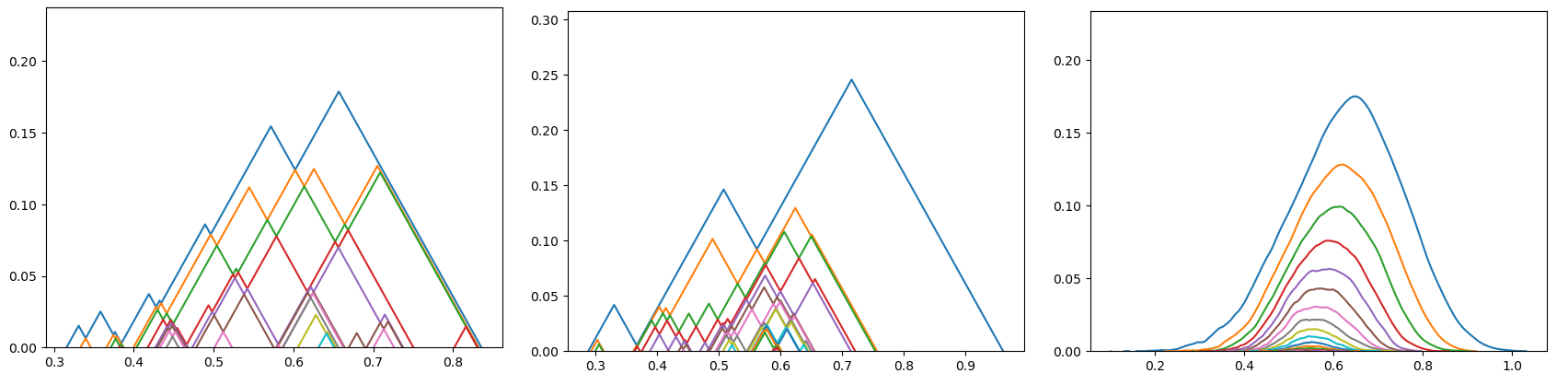}
    \caption{On the left and center: landscape number $4$ and $69$ respectively of the sphere experiment. On the right: the average landscape.}
    \label{landscapemediaesfera}
\end{figure}




\begin{rema}
  From here, when we refer to distance between landscapes, we will be using the \textit{supremum norm} or \textit{infinity norm}, \begin{equation}{\label{supnorm}}
        \norm{\lambda(k,t)-\widetilde{\lambda}(k,t)}=\sup_k\norm{\lambda_k-\widetilde{\lambda}_k}_{\infty}=\sup_k\{\sup_t|\lambda_k(t)-\widetilde{\lambda}_k(t)|\}.
    \end{equation}

    Under this distance, persistence landscapes also have stable behaviour under small perturbations of the original point cloud \cite{bubeniklandscapesstatistical,bubeniklandscapesproperties}.

\end{rema}



\subsection{Computational and data details}

The main code implemented for this work is based on Python language. For persistence diagrams and persistence landscape computations, we have used Ripser \cite{Bauer2021Ripser,ctralie2018ripser} and Persim (\url{https://github.com/scikit-tda/persim}). For some graphical representations we display on this paper we have also used GUDHI (\url{https://gudhi.inria.fr}) \cite{gudhi} and the Landscape Python Package (\url{https://gitlab.com/kfbenjamin/pysistence-landscapes/}). Finally, for the centrality measures computations we have used the Python package Networkx \cite{networkx}.

We have used data sets extracted from Eurocontrol's NEST provided by CRIDA (the Spanish Air Traffic Management R\&D\&I Reference Center). Analyzed data consisted of radar and planned trajectories of all flights with origin or destination at any of the European airports during the period of interest.

The computations performed in this paper were carried out on the servers of the Department of Mathematics at Universidad Autónoma of Madrid and the Department of Quantitative Methods at CUNEF University of Madrid. These resources provided the necessary computational power to handle the complex data analyses and simulations required for our research. 

\section{Application}

As we have explained at the beginning of the paper, we want to use the deviation between the planned and the real trajectory of an airplane as well as its time delay  to obtain a cloud of points of a certain area, airport or, even a country, with the purpose to apply TDA on it.

\subsection{Deviation distance} 

On the way to attach the problem, we have to define a distance to measure the deviation between a planned trajectory and the one the airplane finally did. Each trajectory is given as a set of $4$--dimensional ordered points which have the following configuration: \[
(t=\text{time},L=\text{latitude},l=\text{longitude},a=\text{altitude (km)}).
\]We want to define a distance between two points with that form. For this purpose and to clarify the explanation, suppose we have two points $(t_0,L_0,l_0,a_o)$ and $(t_1,L_1,l_1,a_1)$ where $t_0=t_1$.

\begin{figure}[h]
    \centering
    \includegraphics[width=350pt]{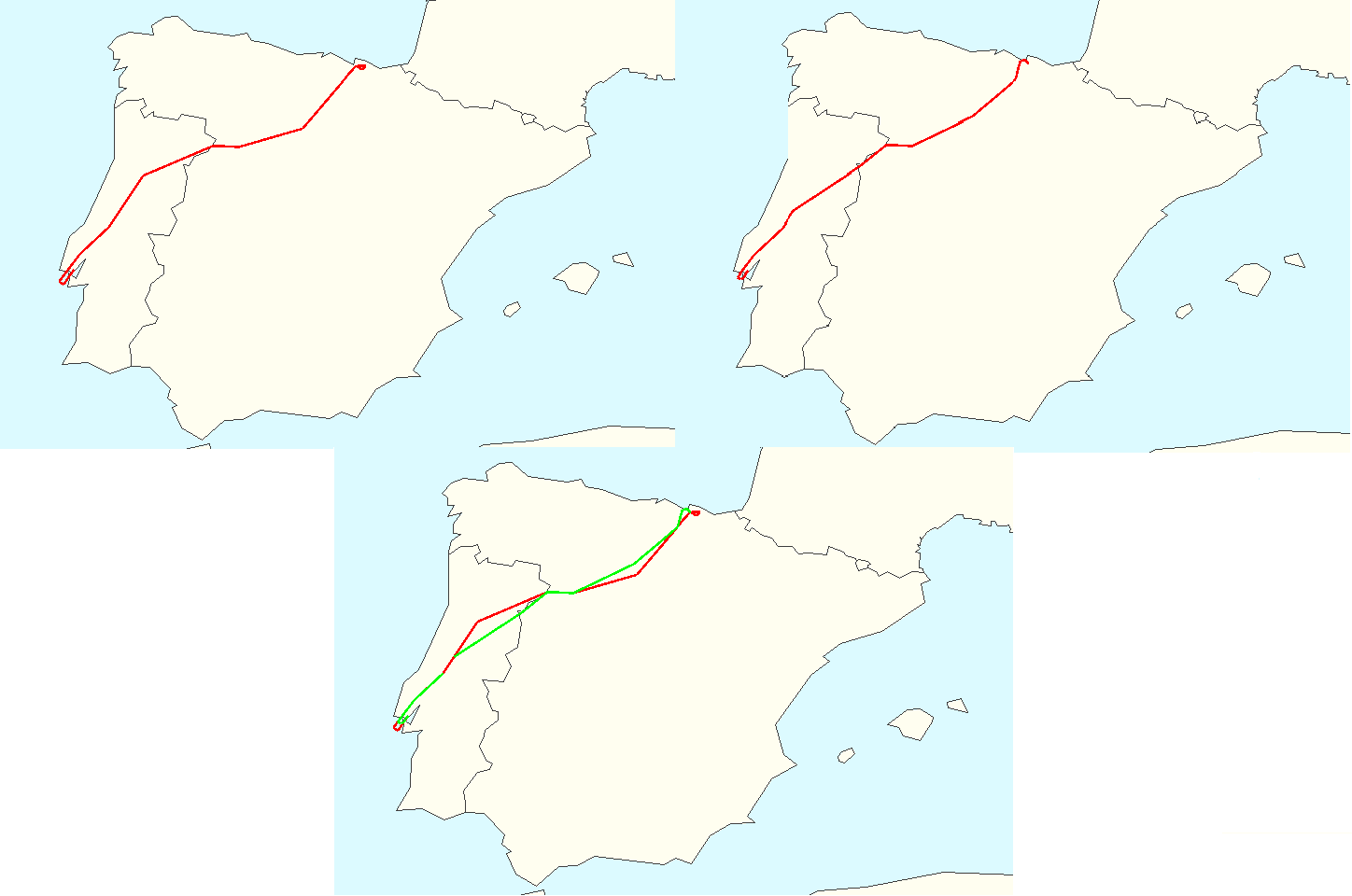}
    \caption{On the top left: The planned trajectory for a Bilbao-Lisbon flight on the 21st July 2019. On the top right: The real trajectory the plane followed. At the bottom: Both trajectories displayed con the same map.}
    \label{santanderlisboa}
\end{figure}

Firstly, we consider the \textit{haversine} distance for the latitude and longitude coordinates.
This calculates the great-circle distance between two points on a sphere based on their longitudes and latitudes. Although the earth is not a perfect round sphere, the haversine distance will be a very good approximation of the real latitude-longitude distance, since we are going to compute mostly distances in Europe. For the altitude coordinate, we simply compute the absolute value of the difference between $a_0$ and $a_1$.

We define a distance between $(t_0,L_0,l_0,a_o)$ and $(t_1,L_1,l_1,a_1)$ where $t_0=t_1=T$ as \begin{equation}{\label{distancia}}
    \dist_{T}((L_0,l_0,a_o),(L_1,l_1,a_1)):=\sqrt{\dist_H((L_0,l_0),(L_1,l_1))^2+|a_0-a_1|^2},
\end{equation}where $\dist_H$ denotes the haversine distance.
This gives a natural, Euclidean distance between the geographical coordinates.

Let $pT$ and $rT$ denote the set of points of a planned and a real trajectory respectively. So, as we have mentioned above, \begin{align*}
    pT=&\{(t_j,L_j,l_j,a_j)\}_{j\in J}\\
    rT=&\{(\widetilde{t}_i,\widetilde{L}_i,\widetilde{l}_i,\widetilde{a}_i)\}_{i\in I},
\end{align*}where $I$ and $J$ are sets of ordered indexes.

\begin{rema}
    There are some considerations about $pT$ and $rT$ we want to make explicit:
    \begin{enumerate}
        \item $(L_{j_0},l_{j_0})=(\widetilde{L}_{i_0},\widetilde{i}_{i_0})$ and $(L_{j_f},l_{j_f})=(\widetilde{L}_{i_f},\widetilde{i}_{i_f})$, where $j_0$, $j_f$, $i_0$ and $i_f$ are the initial and final indexes of each set respectively. This assertion is obvious, as the initial and final locations correspond to the initial and final destinations.
        \item $a_{j_0}=a_{j_f}=a_{i_0}=a_{i_f}=0$.
        \item $\{t_j\}_{j\in J}$ and $\{\widetilde{t}_i\}_{i\in I}$ have normally different cardinality. Even $t_{j_0}\neq\widetilde{t}_{i_0}$ and $t_{j_f}\neq\widetilde{t}_{i_f}$. $\#(\{\widetilde{t}_i\})>\#(\{t_j\})$ as we have much more points in $rT$ than $pT$.
    \end{enumerate}
\end{rema}

Due to the last point of the Remark and the definition \eqref{distancia} of our distance, for computing the distance between $pT$ and $rT$ we have, firstly, to redefine those sets so that they share the same temporal marks. Our choice for the new temporal sequence will be just taking the union of the temporal sequences for the planned and the real flight: 
\[
\{\widehat{t}_k\}_{k\in K}=\{t_j\}_{j\in J}\cup\{\widetilde{t}_i\}_{i\in I},
\]where $\widehat{t}_{k_0}=\min\{t_{j_0},\widetilde{t}_{i_0}\}$ and $\widehat{t}_{k_f}=\max\{{t_{j_f},\widetilde{t}_{i_f}}\}$, where $k_0$ denotes the first index (we will  identify for clarity $k_0=0$) of $K$ and $k_f$ the last one.

Now, without loss of generality, suppose that $\widehat{t}_{k_l}\notin\{t_j\}_{j\in J}$ and $\widehat{t}_{k_{l-1}},\widehat{t}_{k_{l+1}}\in\{t_j\}_{j\in J}$. We need to create a new point $p=(\widehat{t}_{k_l},L_{k_l},l_{k_l},a_{k_l})$. For that purpose, we have chosen $p$ on the segment between $(\widehat{t}_{k_{l-1}},L_{k_{l-1}},l_{k_{l-1}},a_{k_{l-1}})$ and $(\widehat{t}_{k_{j+1}},L_{k_{l+1}},l_{k_{l+1}},a_{k_{l+1}})$ according to the nature of our data . Furthermore, if $\widehat{t}_{k_{l-1}}\in\{t_j\}_{j\in J}$ but $\widehat{t}_{k_{l+1}}\notin\{t_j\}_{j\in J}$, we pick the first $\widehat{t}_{k_p}\in\{t_j\}_{j\in J}$ such that $k_p>k_l$ and interpolate between $k_p$ and $k_l$ to obtain the required information.

In addition,  suppose $t_{j_0}>\widehat{t}_0$ and $t_{j_f}<\widehat{t}_{k_f}$. Let $k_{j_0}$, $k_{j_f}\in K$ such that $\widehat{t}_{k_{j_0}}=t_{j_0}$ and $\widehat{t}_{k_{j_f}}=t_{j_f}$; then $(\widehat{t}_{k_l},L_{k_l},l_{k_l},a_{k_l})=(\widehat{t}_{k_l},L_{j_0},l_{j_0},a_{j_0})$ and $(\widehat{t}_{k_p},L_{k_p},l_{k_p},a_{k_p})=(\widehat{t}_{k_p},L_{j_f},l_{j_f},a_{j_f})$, for $k_l<k_{j_0}$ and $k_p>k_{j_f}$.

So, after this construction, we have redefined our trajectories and settled $pT$ and $rT$ on the same temporal sequence. Now, it has sense to apply the distance defined in \eqref{distancia} to every pair or points of $pT\times rT$. We can finally define the distance which measures the deviation (in kilometres) between a planned and a real trajectory in the following manner \begin{equation}{\label{distanciasuma}}
    \dist_{\text{Dev}}(pT,rT):=\sum_{l=0}^{k_f}\dist_{\widehat{t}_l}((L_l,l_l,a_l),(\widetilde{L}_l,\widetilde{l}_l,\widetilde{a}_l)).
\end{equation}

\subsection{Point cloud}

As shown in Section \ref{methods}, in order to apply TDA techniques to a given problem, we require a point cloud to work with. To better illustrate this point, let us consider a particular airport that we wish to focus on, and investigate how deviations and delays affect its performance. Suppose that this airport has ten flights per day. As discussed in the previous subsection, we compute the distances between the planned and actual trajectories of each flight, resulting in a list of ten distances denoted as $\{d_1,\dots,d_{10}\}$. We can define our point cloud as follows:\[
\mathcal{V}=\{(\widetilde{d}_i,r_i)\}_{i\in[1,10]}\subset\RR^2,
\] 
where $\widetilde{d}_i=(-1)^pd_i$ with $p=0$ if the flight arrived late or $p=1$ if the flight arrived on time or sooner, and $r_i$ corresponds to the quantity $\widetilde{t}_{i_f}-t_{j_f}$, i.e., the difference between the time the flight landed and the expected one. If the flight was late, we will have $r_i>0$, and on the other hand, $r_i\leq0$ if the flight was on time or sooner. The points will be settled in the first and third quadrant of $\RR^2$ as $\mathcal{V}$ only has points with both positive or bot negative coordinates.

\begin{figure}[h]
    \centering
    \includegraphics[width=400pt]{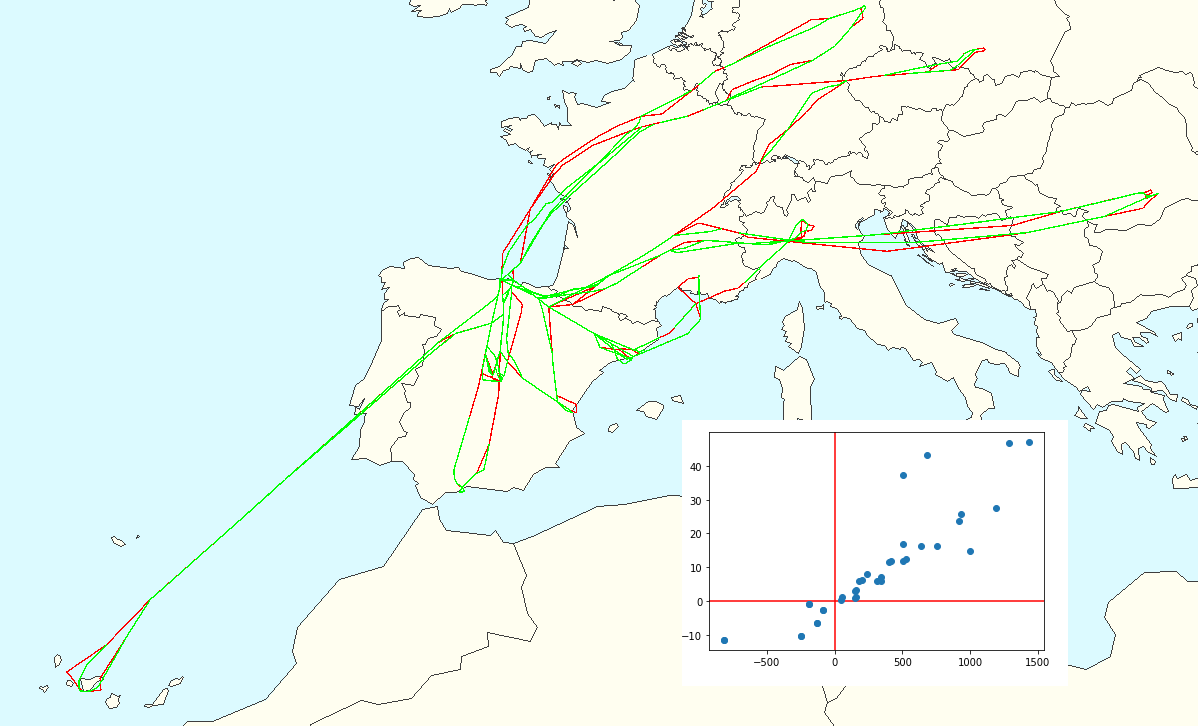}
    \caption{Point cloud of Santander's airport on the 29th June 2019.}
    \label{pointcloudsantander}
\end{figure}

In this paper, we are focusing on airports in order to obtaining the point clouds. But, as this construction is generic, we can choose more than one airport, a certain region or even a country.

\subsection{Persistent homology}

Once we have built our point cloud we can apply TDA on it. In this case, we will work with the $0$--homology as we want to extract the information of how the connected components in the Vietoris-Rips of the filtration of our point cloud $\mathcal{V}$ birth and die. With this information we recover how the points are disseminated on $\RR^2$. Intuitively, if one flight is highly delayed, this point would be far from the origin $(0,0)$ and that phenomena would be detected by the persistent homology analysis. 

\begin{figure}[h]
    \centering
    \includegraphics[width=250pt]{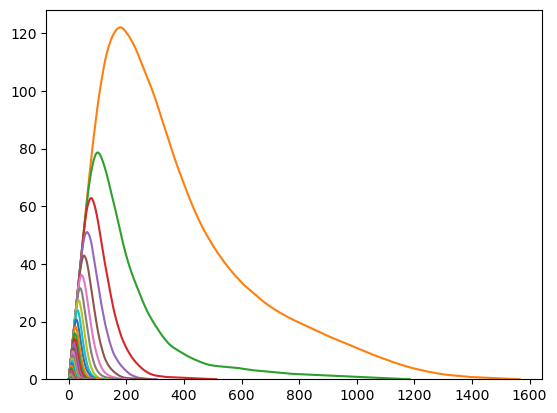}
    \caption{Average Persistence Landscape of Santander's airport in the Summer Season of 2018.}
    \label{PLaverageSantanderSS2018}
\end{figure}

For the explanation of the process we have developed, we are going to illustrate it with a particular example: Santander's airport (LEXJ) on the Summer Season (23rd March to 27th October) 2018. Although the data of 2018 can be considered old, it is not the case for the sake of this study, as current flight patterns are not significative different from those of 2018. For every day during this period of time, we are going to compute its corresponding point cloud $\mathcal{V}$ and extract from there the $PD_0(\mathcal{V})$. As the Summer Season of 2018 has $217$ days, after this first step we will obtain a set $V=\{PD_0(\mathcal{V}_j)\}$ with $217$ persistence diagrams; one for each day. After this, we compute the persistence landscape of each persistence diagram and compute their average persistence landscape (Figure \ref{PLaverageSantanderSS2018}).

\begin{figure}[!h]
    \centering
    \includegraphics[width=\columnwidth]{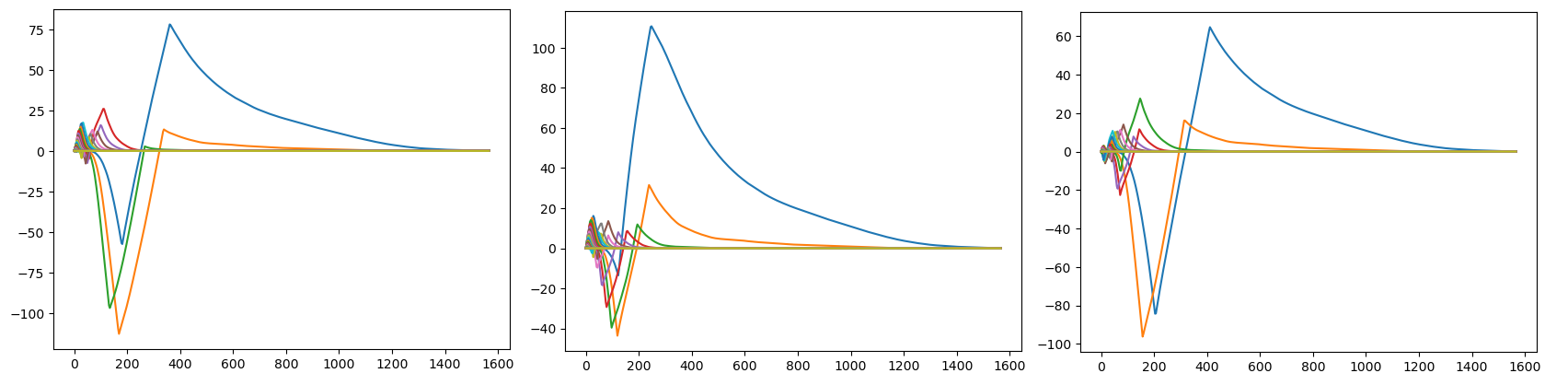}
    \caption{The nine point clouds with persistence landscapes with biggest distance with respect to the average persistence landscape of Santander's airport in the Summer Season of 2018.}
    \label{landscape diferencia}
\end{figure}

\begin{figure}[!h]
    \centering
    \includegraphics[width=\columnwidth]{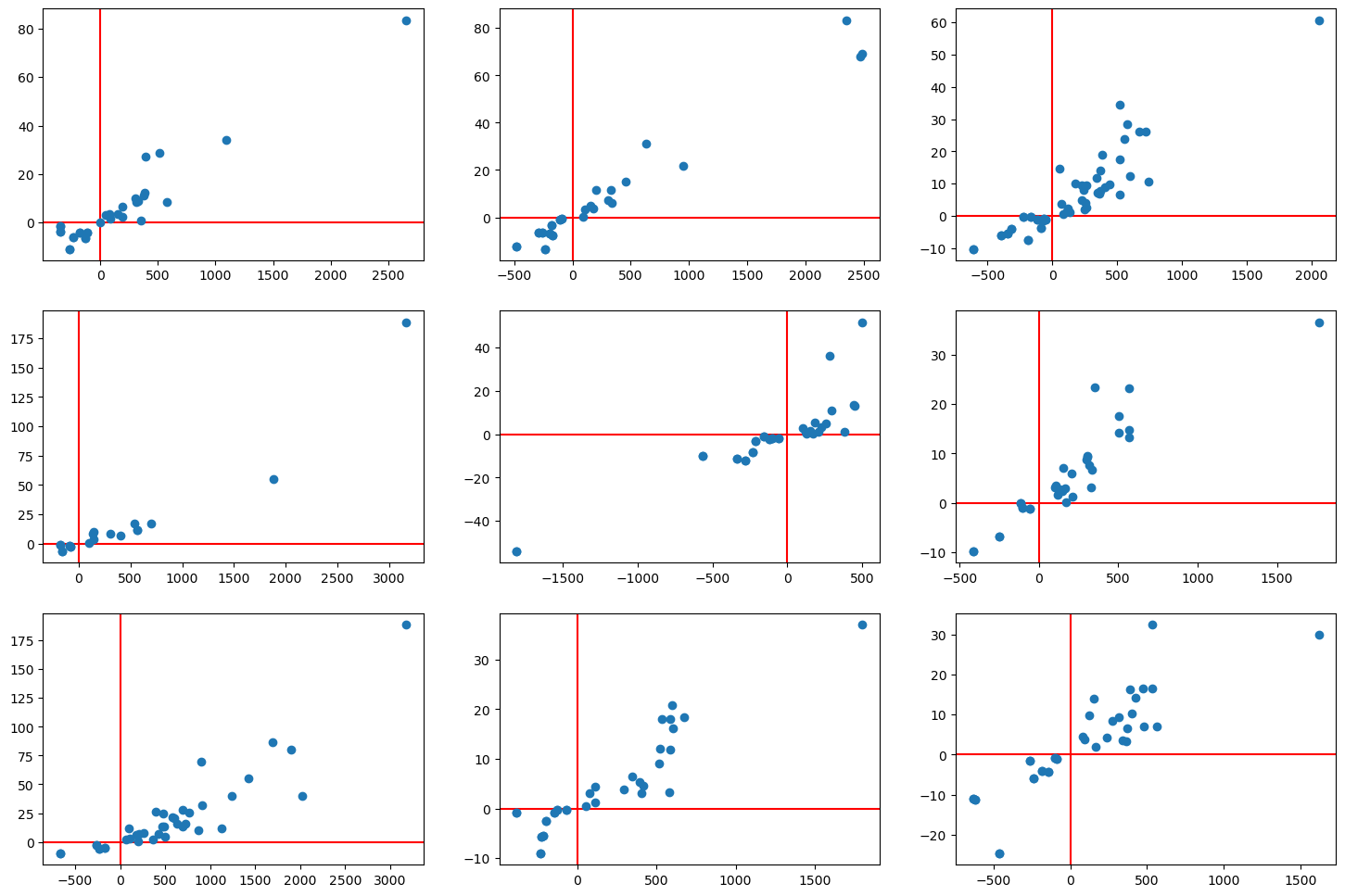}
    \caption{The nine point clouds with persistence landscapes with biggest distance with respect to the average persistence landscape of Santander's airport in the Summer Season of 2018.}
    \label{9pointcloud}
\end{figure}

With this average persistence landscape we obtain a particular image of every airport depending on the distribution of its daily flights via trajectory deviation and delay. One of the advantages of landscapes is that we can now compute the distance with the supremum norm we illustrated on \eqref{supnorm} of each day with respect to the average (Figure \ref{landscape diferencia}). After performing these computations, the point clouds associated to the days with biggest distance from the average have clearly isolated points (Figure \ref{9pointcloud}).

To summarize, we present the pseudocode of the algorithm explained in this section (Algorithm \ref{alg:footprint}) to clarify the process for the rest of the paper, before presenting the results.

\begin{algorithm}
\caption{Algorithm for the footprint of an airport}
\label{alg:footprint}
\begin{algorithmic}
\Require $N\geq 0$ and \textbf{airport} \Comment{\textit{$N$ is the number of days computed in certain airport}}
\State Create an empty list \textbf{PD} for the persistence diagrams
\For{day in $N$}
    \State Create lists with planned \textbf{PT} and real trajectories \textbf{RT} of that \textbf{day}
    \State Create an empty list \textbf{PCD} for the point cloud of the day
    \For{flight in \textbf{PT} and \textbf{RT}}
        \State Compute (deviation distance, delay) and add the point to \textbf{PCD} 
    \EndFor
    \State Compute the distance matrix \textbf{DMD} of \textbf{PCD}
    \State Compute the persistence diagram \textbf{PDD} of \textbf{DMD} and add it to \textbf{PD} \Comment{\textit{For this step, we use Ripser package}}
\EndFor
\State Create an empty list \textbf{PL} for the presistence landscapes
\For{\textbf{pd} in \textbf{PD}}
    \State Compute de persistence landscape \textbf{PLD} of \textbf{pd} and add it to \textbf{PL} \Comment{\textit{For this step, we use Persim package}}
\EndFor
\State Finally, we computed de average persistence landscape \textbf{APL} (footprint) of \textbf{PL} \Comment{\textit{For this step, we also use Persim package}}

\end{algorithmic}
\end{algorithm}

\section{Discussion and results}

The purpose of this paper is to introduce the concepts of Topological Data Analysis (TDA) and explore their application in the context of Air Traffic Management (ATM) via a proof of concept presented in this section. We aim to demonstrate the efficacy of TDA through an analysis of real-world data, specifically the Spanish network of airports during the Summer Season of 2018, as classified by AENA (a Spanish public company responsible for overseeing airports of general interest in Spain).

In order to facilitate our analysis, we have leveraged AENA's airport classification system for 2018, which is divided into five distinct categories. These categories include Group 3 (airports with less than half a million passengers per year), Group 2 (airports with between half a million and two million passengers per year), Group 1 (airports with more than two million passengers per year), the Canary Group (comprised of airports located in the Canary Islands), and a special group (which includes Madrid-Barajas, Barcelona-El Prat, and Palma's airport).


\begin{rema}
    The upper bound imposed on Group 3 necessitated the inclusion of airports of varying types. In light of this, we have partitioned Group 3 into four distinct subgroups, namely general aviation airports, air bases open to civilian traffic, airports with low traffic, and helipads (which will be excluded from our analysis).

    Similarly, the Canary Group exhibits a significant variation in the size of its airports, ranging from La Gomera with two to four flights per day to Gran Canaria with over one hundred flights per day. Despite this disparity, we have chosen to maintain all airports in the same group for our analysis.
\end{rema}

From now on, we will denote each airport by its ICAO (International Civil Aviation Organization) code. Here is the classification given by AENA in 2018:

\begin{itemize}
    \item \textit{Group 3}: \begin{itemize}
        \item \textit{General aviation airports}: Madrid-C. Vientos (LECU) and Sabadell (LELL).
        \item \textit{Air bases open to civilian traffic}: Albacete (LEAB), Badajoz (LEBZ), León (LELN), Salamanca (LESA), Son Bonet (LESB) and Valladolid (LEVD).
        \item \textit{Airports with low traffic}: Melilla (GEML), Burgos-Villafría (LEBG), Córdoba (LEBA), Girona (LEGE), Huesca-Pirineos (LEHC), Logroño-Agoncillo (LERJ),  Pamplona (LEPP), San Sebastián (LESO) and Vitoria (LEVT).
    \end{itemize}
    \item \textit{Group 2}: A Coruña (LECO), Almería (LEAM), Asturias (LEAS), Federico García Lorca Granada-Jaén (LEGR), Jerez (LEJR), Reus (LERS), Santander (LEXJ), Santiago (LEST), Vigo (LEVX) and \textit{Zaragoza (LEZG)}.
    \item \textit{Group 1}: Alicante Elche-Miguel Hernández (LEAL), Bilbao (LEBB), Ibiza (LEIB), Málaga-Costa del Sol (LEMG), Menorca (LEMH), Sevilla (LEZL) and Valencia (LEVC).
    \item \textit{Canary Group}: Fuerteventura (GCFV), La Gomera (GCGM), El Hierro (GCHI),  La Palma (GCLA),  Gran Canaria (GCLP), César Manrique Lanzarote (GCRR),  Tenerife South (GCTS) and Tenerife North (GCXO).
    \item \textit{Special Group}: Adolfo Suárez Madrid-Barjas (LEMD), Josep Tarradellas Barcelona-El Prat (LEBL) and Palma Mallorca (LEPA).
\end{itemize}

\begin{rema}
    We decided to highlight Zaragoza's airport (in italic) due to its uniqueness. It is the only airport in Group 2 that does not facilitate passenger aviation. We will see that this phenomena is detected by our analysis.
\end{rema}

\begin{figure}[!h]
    \centering
    \includegraphics[width=\columnwidth]{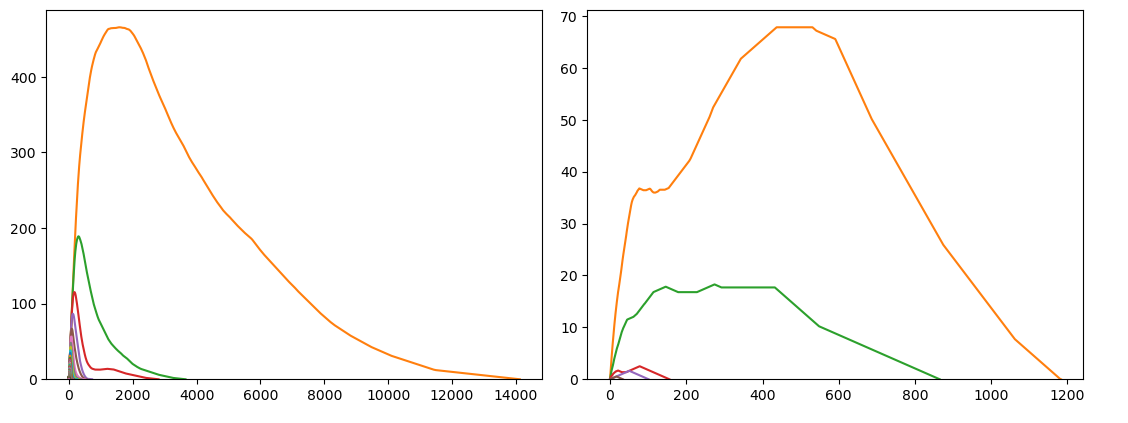}
    \caption{On the left: Average Persistence landscape of Adolfo Suárez Madrid-Barjas' airport on the Summer Season of 2018. On the right: Average Persistence landscape of Huesca-Pirineos' airport on the Summer Season of 2018. Check the difference in magnitude of the axes}
    \label{PLhuescamadrid}
\end{figure}

First, we compute the average persistence landscape of every airport corresponding to the Summer Season of 2018 obtaining a collection of average persistence landscapes $\mathcal{L}=\{\overline{\lambda}(k,t,A)\}_{A\in\mathcal{A}}$, where $\mathcal{A}$ represents the list of all of the Spanish airports. Notably, there are a wide variety of landscapes in $\mathcal{L}$. This type of analysis benefits the number of flights each airport has. So, airports with big differences in the number of passengers per year would have highly different landscapes (see Figure \ref{PLhuescamadrid} for a comparison between Adolfo Suárez Madrid-Barajas' and Huesca-Pirineos paying attention the grids of the pictures).


Second, we construct a distance matrix $\mathcal{D}=(d_{ij})$, where each $d_{ij}$ is the supremum distance (explained on \eqref{supnorm}) between the persistence average landscape of the airport $i$ and the airport $j$. We have split $\mathcal{D}$ in different tables in order to clearly present the relevant information. We present the most interesting considerations regarding them. For the interested reader, the rest of the tables are displayed on Appendix \ref{tablas}:

\begin{itemize}
    \item On Table \ref{tablagrupo3BT} we have highlighted four airports: GEML (red), LEGE (blue), LEPP (blue) and LEVT (blue). On Table \ref{tablagrupo23BT} we have higlighted three airports: LEGE, LEPP and LEVT (all in blue).

        \begin{table}[h]
    \centering
\resizebox{\columnwidth}{!}{%
\begin{tabular}{l|ccccccccc}
{} & GEML & LEBA & LEBG & LEGE & LEHC & LEPP & LERJ & LESO & LEVT \\
\midrule
GEML & 0.00 & 93.78 & 69.05 & 117.80 & 67.61 & 105.05 & 62.34 & 66.97 & 109.09 \\
LEBA & 93.78 & 0.00 & 30.87 & 66.68 & 65.14 & 42.23 & 31.76 & 32.52 & 61.30 \\
LEBG & 69.05 & 30.87 & 0.00 & 77.77 & 43.86 & 53.81 & 30.38 & 40.65 & 72.06 \\
LEGE & 117.80 & 66.68 & 77.77 & 0.00 & 86.07 & 35.63 & 67.72 & 57.29 & 23.91 \\
LEHC & 67.61 & 65.14 & 43.86 & 86.07 & 0.00 & 70.16 & 34.20 & 45.35 & 80.90 \\
LEPP & 105.05 & 42.23 & 53.81 & 35.63 & 70.16 & 0.00 & 46.40 & 45.63 & 25.94 \\
LERJ & 62.34 & 31.76 & 30.38 & 67.72 & 34.20 & 46.40 & 0.00 & 36.44 & 63.69 \\
LESO & 66.97 & 32.52 & 40.65 & 57.29 & 45.35 & 45.63 & 36.44 & 0.00 & 45.05 \\
LEVT & \boxit{0.45in}{1.8in}109.09 & 61.30 & 72.06 & \boxittwo{0.45in}{1.8in} \ 23.91 & 80.90 & \boxittwo{0.45in}{1.8in} \ 25.94 & 63.69 & 45.05 & \boxittwo{0.45in}{1.8in} \ \ 0.00 \\
\end{tabular}}
    \caption{Distance matrix corresponding to the Group 3 - \textit{Airports with low traffic}}
    \label{tablagrupo3BT}
\end{table}

\begin{table}[h]
    \centering
\resizebox{\columnwidth}{!}{%
\begin{tabular}{l|ccccccccc}
{} & GEML & LEBA & LEBG & LEGE & LEHC & LEPP & LERJ & LESO & LEVT \\
\midrule
LEAM & 116.56 & 68.71 & 79.58 & 35.88 & 87.21 & 37.49 & 70.24 & 53.22 & 14.89 \\
LEAS & 101.64 & 59.92 & 68.69 & 44.71 & 75.33 & 26.39 & 61.26 & 39.34 & 22.59 \\
LECO & 94.07 & 47.88 & 55.34 & 64.26 & 67.45 & 41.22 & 49.55 & 27.13 & 40.81 \\
LEGR & 102.80 & 58.45 & 66.47 & 56.37 & 75.31 & 34.86 & 59.47 & 35.85 & 32.92 \\
LEJR & 120.01 & 63.52 & 74.16 & 24.02 & 87.82 & 31.04 & 64.35 & 57.16 & 13.78 \\
LERS & 116.04 & 69.39 & 80.49 & 36.89 & 86.45 & 36.90 & 70.73 & 55.17 & 15.85 \\
LEST & 113.49 & 60.70 & 68.37 & 27.82 & 81.63 & 29.77 & 60.70 & 49.35 & 11.10 \\
LEVX & 101.06 & 54.88 & 61.78 & 45.00 & 70.76 & 24.40 & 55.77 & 36.73 & 21.64 \\
LEXJ & 113.84 & 61.08 & 68.05 & 44.10 & 83.30 & 29.01 & 60.72 & 48.26 & 20.35 \\
LEZG & 242.82 & 202.76 & 207.87 & \boxittwo{0.46in}{2in}157.17 & 181.63 & \boxittwo{0.46in}{2in}185.35 & 188.05 & 199.45 & \boxittwo{0.48in}{2in}174.77 \\
\end{tabular}}
    \caption{Distance matrix corresponding to the Group 2 (rows) and Group 3 - \textit{Airports with low traffic} (columns) where we have remarked two columns: Girona's and Vitorias's airport respectively}
    \label{tablagrupo23BT}
\end{table}

Melilla's airport (GEML) seems very different from the ones in its group and the ones in Group 2. It seems  that its location affects it in comparison to the rest of the airports in the Iberica's Peninsula. 

The airports highlighted in blue seems to work better with Group 2 than with Group 3 - \textit{Airports with low traffic} in the Summer Season 2018. A few years later, Girona's airport (LEGE) was indeed included in Group 2. 

\item Table \ref{tablagrupo2} represents the symmetric part of $\mathcal{D}$ corresponding to Group 2.

\begin{table}[h]
\centering
\resizebox{\columnwidth}{!}{%
\begin{tabular}{l|cccccccccc}
{} & LEAM & LEAS & LECO & LEGR & LEJR & LERS & LEST & LEVX & LEXJ & LEZG \\
\midrule
LEAM & 0.00 & 20.27 & 34.82 & 26.39 & 21.01 & 15.06 & 18.62 & 29.54 & 15.59 & 178.62 \\
LEAS & 20.27 & 0.00 & 22.89 & 15.39 & 23.99 & 19.35 & 24.05 & 11.82 & 13.71 & 195.42 \\
LECO & 34.82 & 22.89 & 0.00 & 14.46 & 42.92 & 36.61 & 37.09 & 21.62 & 24.81 & 212.38 \\
LEGR & 26.39 & 15.39 & 14.46 & 0.00 & 36.74 & 30.78 & 29.65 & 13.62 & 20.32 & 204.78 \\
LEJR & 21.01 & 23.99 & 42.92 & 36.74 & 0.00 & 13.80 & 16.74 & 29.64 & 20.73 & 175.38 \\
LERS & 15.06 & 19.35 & 36.61 & 30.78 & 13.80 & 0.00 & 19.84 & 31.14 & 16.62 & 188.94 \\
LEST & 18.62 & 24.05 & 37.09 & 29.65 & 16.74 & 19.84 & 0.00 & 19.25 & 16.69 & 177.77 \\
LEVX & 29.54 & 11.82 & 21.62 & 13.62 & 29.64 & 31.14 & 19.25 & 0.00 & \boxitthree{0.45in}{0.30in}16.14 & 191.58 \\
LEXJ & 15.59 & \boxitthree{0.45in}{0.30in}13.71 & 24.81 & 20.32 & 20.73 & 16.62 & 16.69 & 16.14 & 0.00 & 193.29 \\
LEZG & 178.62 & 195.42 & 212.38 & 204.78 & 175.38 & 188.94 & 177.77 & 191.58 & 193.29 & \boxit{0.46in}{2in} \ \ 0.00 \\
\end{tabular}}
\caption{Distance matrix corresponding to the Group 2 airports where we have remarked the Zaragoza's airport column and some close distances} 
\label{tablagrupo2}
\end{table}

As we mentioned before, Zaragoza's airport (LEZG) is the only airport in Group 2 that does not facilitate passenger aviation. That pathological behaviour is clearly detected by the TDA analysis. 

Furthermore, we have observed that airports located in close proximity tend to have similar geometric properties, as it is the case with Asturias (LEAS), Santander (LEXJ), Santiago (LEST), and Vigo (LEVX). This phenomenon suggests that the point cloud we create based on deviation of trajectories and delays, encodes information regarding geographical location of the airports and origins and destination of its flights. TDA analysis is capable of capturing that as we can check in the distances of airports with those similar characteristics.


\item On Table \ref{tablagrupo12} we show the differences between Group 2 and Group 1.

\begin{table}[h]
    \centering
\resizebox{\columnwidth}{!}{%
\begin{tabular}{l|cccccccccc}
{} & LEAM & LEAS & LECO & LEGR & \boxitthree{0.42in}{0.14in}LEJR & LERS & \boxit{0.42in}{0.14in}LEST & LEVX & LEXJ & LEZG \\
\midrule
LEAL & 65.50 & 75.44 & 94.91 & 87.04 & 54.71 & 68.10 & 58.41 & 75.54 & 74.80 & 128.03 \\
LEBB & 29.20 & 41.97 & 59.93 & 52.40 & 22.80 & 35.39 & \boxit{0.4in}{0.11in}24.14 & 39.74 & 40.26 & 154.06 \\
LEIB & 31.66 & 48.49 & 65.37 & 57.83 & 29.56 & 42.30 & 31.02 & 44.62 & 46.30 & 147.48 \\
LEMG & 101.62 & 113.41 & 132.27 & 124.55 & 92.45 & 105.89 & 95.73 & 112.71 & 112.38 & 98.78 \\
\boxittwo{6.1in}{0.11in}LEMH & 22.24 & 13.10 & 23.64 & 18.86 & 22.87 & 16.99 & \boxit{0.4in}{0.11in}24.12 & 18.83 & 9.37 & 197.85 \\
LEVC & 58.15 & 64.32 & 83.97 & 76.14 & 43.39 & 56.81 & 47.25 & 64.78 & 63.47 & 140.87 \\
\boxitthree{0.41in}{0.11in}LEZL & 28.20 & 33.84 & 52.66 & 44.80 & \boxitthree{0.4in}{0.11in}15.41 & 27.62 & \boxit{0.4in}{0.11in}22.85 & 33.74 & 31.92 & \boxittwo{0.5in}{1.47in}161.61 \\
\end{tabular}}
    \caption{Distance matrix corresponding to the Group 2 (columns) and Group 1 (rows) with different marks}
    \label{tablagrupo12}
\end{table}

We can see that Menorca's airport (LEMH) is close to the airports of Group 2. Moreover, its distance to Santander's airport (LEXJ) is one of the closest in $\mathcal{D}$. 

Nowadays, Santiago's airport (LEST) is in Group 1. We can check how in 2018 this aiport was close to some airports (Bilbao (LEBB), Menorca (LEMH), Sevila (LEZL)) of Group 1.

Finally, we want to highlight the similarity between Sevilla's (LEZL) and Jerez's (LEJR) airport (here it is another example of two airports that are very close geographically) and how Zaragoza's aiport (LEZG) due to its uniqueness it is still very different from airports of Group 1.

Accordingly to the information displayed above, on table \ref{tablagrupo1} we can see the symmetric part of $\mathcal{D}$ related to Group 1, and how Menorca's airport (LEMH) has bigger distances with its own group than with Group 2.

\begin{table}[h]
    \centering
\resizebox{380pt}{!}{%
\begin{tabular}{l|ccccccc}
{} & LEAL & LEBB & LEIB & LEMG & LEMH & LEVC & LEZL \\
\midrule
LEAL & 0.00 & 37.29 & 37.18 & 38.29 & 77.34 & 17.64 & 42.97 \\
LEBB & 37.29 & 0.00 & 19.61 & 73.43 & 44.16 & 29.01 & 26.18 \\
LEIB & 37.18 & 19.61 & 0.00 & 73.16 & 51.19 & 29.37 & 21.39 \\
LEMG & 38.29 & 73.43 & 73.16 & 0.00 & 115.33 & 52.71 & 80.47 \\
LEMH & 77.34 & 44.16 & 51.19 & 115.33 & 0.00 & 66.20 & 36.46 \\
LEVC & 17.64 & 29.01 & 29.37 & 52.71 & 66.20 & 0.00 & 31.66 \\
LEZL & 42.97 & 26.18 & 21.39 & 80.47 & \boxit{0.42in}{1.45in}36.46 & 31.66 & 0.00 \\
\end{tabular}}
    \caption{Distance matrix corresponding to the Group 1 with Menorca's airport highlighted}
    \label{tablagrupo1}
\end{table}

\item Canary Group is interesting due to two important factors: its geographical location faraway from the rest airports of Spain and the diversity on the size of the airports in the group. 

\begin{table}[h]
    \centering
\resizebox{\columnwidth}{!}{%
\begin{tabular}{l|cccccccc}
{} & GCFV & GCGM & GCHI & GCLA & GCLP & GCRR & GCTS & GCXO \\
\midrule
GCFV & 0.00 & 153.65 & 156.43 & 31.55 & 37.29 & 27.18 & 78.66 & 80.91 \\
GCGM & 153.65 & 0.00 & \boxittwo{0.40in}{0.1in}11.09 & 163.19 & 144.32 & 141.84 & 183.45 & 102.31 \\
GCHI & 156.43 & \boxittwo{0.43in}{0.1in}11.09 & 0.00 & 166.15 & 147.30 & 144.71 & 185.81 & 105.02 \\
GCLA & 31.55 & 163.19 & 166.15 & 0.00 & 63.54 & 23.46 & 110.20 & 72.97 \\
GCLP & 37.29 & 144.32 & 147.30 & 63.54 & 0.00 & 53.73 & 64.41 & 69.07 \\
GCRR & 27.18 & 141.84 & 144.71 & 23.46 & 53.73 & 0.00 & 88.74 & 62.67 \\
GCTS & 78.66 & 183.45 & 185.81 & 110.20 & 64.41 & 88.74 & 0.00 & 126.77 \\
GCXO & 80.91 & \boxit{1.03in}{1.63in}102.31 & 105.02 & 72.97 & 69.07 & 62.67 & 126.77 & 0.00 \\
\end{tabular}}
    \caption{Distance matrix corresponding to the Canary Group with different marks}
    \label{tablagrupoc}
\end{table}

La Gomera's and El Hierro's airport are the smallest of that group (for example, La Gomera's airport has four flights per day and El Hierro approximately ten). Its mutual distance is very small, as we can see on Table \ref{tablagrupoc}, in comparison with the rest of the Canary airports.

 It is very interesting to check how the Canary Group behaves with the rest of Spanish airport net. It seems that its geographical location and the length of its flights (on longer flights the deviation between trajectories increases) are an important factor to characterize this group. 

 \item On Table \ref{tablagrupoe} we displayed the symmetric part of $\mathcal{D}$ corresponding to the Special Group.

 \begin{table}[h]
    \centering
\begin{tabular}{l|ccc}
{} &        LEBL &        LEMD &        LEPA \\
\midrule
LEBL &    0.00 &  251.13 &  337.09 \\
LEMD &  251.13 &    0.00 &  463.99 \\
LEPA &  337.09 &  463.99 &    0.00 \\
\end{tabular}
    \caption{Distance matrix corresponding to the Special Group}
    \label{tablagrupoe}
\end{table}

Due to the uniqueness of this airports, the distance between them are huge and disparate. The information arose from this analysis suggests that Palma's airport (LEPA) is more similar to the airports from the existent groups than to Madrid or Barcelona. This can be checked in the tables on Appendix \ref{tablas}, but for the clarity of the remark, we displayed on Table \ref{tablagrupo1e} the part of $\mathcal{D}$ corresponding to the Special Group with Group 1.

\begin{table}[h]
    \centering
\begin{tabular}{l|ccc}
{} &        LEBL &        LEMD &       LEPA \\
\midrule
LEAL &  306.89 &  448.50 &  36.36 \\
LEBB &  324.83 &  460.19 &  21.29 \\
LEIB &  317.48 &  440.60 &  23.69 \\
LEMG &  281.08 &  438.53 &  71.78 \\
LEMH &  360.30 &  465.55 &  47.73 \\
LEST &  336.21 &  459.76 &  32.27 \\
LEVC &  313.41 &  453.90 &  23.77 \\
LEZL &  329.62 &  461.41 &  \boxit{0.44in}{1.63in}17.11 \\
\end{tabular}
    \caption{Distance matrix corresponding to the Special Group (columns) and the Group 1 (rows) with Palma's airport (LEPA) highlighted}
    \label{tablagrupo1e}
\end{table}

\end{itemize}

This detailed study of the distance matrix $\mathcal{D}$ has clearly illustrated how much intricate information TDA is able to extract from an apparently very easy point cloud based on deviation of trajectories and delays. We found that, apart for the number of passengers per year in each airport on which the Spanish classification is now constructed on, geographical properties, origins and destinations and typology of the airport also influence and appears on the classification that TDA shows.

We also studied the statistical significance of the distances obtained in $\mathcal{D}$. To do this, and following the statistical tests performed by Bubenik and Dłotko in \cite{bubeniklandscapesstatistical,bubeniklandscapestoolbox}, we decided to compute a \textit{block permutation test} \cite{lahiri2013resampling} between each airport. Instead of permuting all persistence landscapes of airport $i$ with those of airport $j$, we permuted blocks of consecutive landscapes of airport $i$ with blocks of consecutive landscapes of airport $j$, considering the dependencies between consecutive days in air traffic. Due to computational limits, we calculated 400 permutations for each pair of airports.

Although the complete results of the $p$-values can be found in Appendix \ref{app:p-value}, we present some interesting conclusions extracted from them. As we previously pointed out, Palma's airport (LEPA) would fit better in other groups such as Group 1 or 2 (either some big $p$-values with airports of Group 3), at least if we base our classification on trajectories and delays. Additionally, the square parts of Table \ref{table:p-valuegroups12} regarding Group 2 and Group 1 respectively, reinforce our conclusions about the similarities among them in terms of air traffic, and also between groups. Finally, the variability in capacity and number of flights per day in airports of Group 3 (some of them, like LEHC, had only 31 flights in the entire 2018 Summer Season) is reflected in their $p$-values. They also follow the pattern of having higher $p$-values with airports of the same group, but some anomalies are detected when comparing them with Group 1. We assume that the difference in capacity and number of flights per day could be the reason for these discrepancies.

\begin{figure}[!h]
    \centering
    \includegraphics[width=220pt]{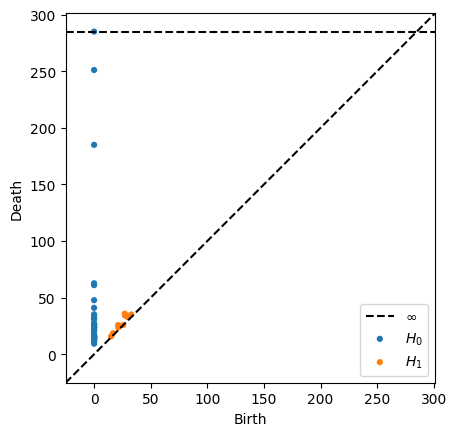}
    \caption{0- and 1-persistence diagrams of the net of Spanish airports in the Summer Season of 2018.}
    \label{fig:PD spanish airport net}
\end{figure}

A final step regarding matrix $\mathcal{D}$ would be trying to replicate a point cloud in some Euclidean space with a distance matrix $\widetilde{\mathcal{D}}$ similar to $\mathcal{D}$. There are a lot of different and successful methods such as Isomap \cite{isomap}. Instead of doing that, we consider that applying TDA as well to $\mathcal{D}$ (Ripser allows a distance matrix as input for computations) is suitable for this kind of experiments. As we explained, TDA tries to deal with high dimensionality and interconnectedness, so a natural question will be how the points that realize $\mathcal{D}$ as distance matrix behave. In future works, we hope to explore and present this part of the research, trying to provide some classification to the whole European airport net. Meanwhile, on Figure \ref{fig:PD spanish airport net} the reader can see the 0- and 1-persistence diagram produced by the Spanish airport network during the Summer Season of 2018.

Based on the conclusions extracted from $\mathcal{D}$ and the $p$-value computations, we see that improving the data used to classify airports (the Spanish Airport Network classification currently relies on the number of passengers per year) would provide a more complex and detailed picture of the network. This, in turn, would help in taking similar actions at airports of the same class and reducing redundant operations. We firmly believe that adding additional information, such as the order of flights per day or even their capacities, would enrich the experiment we have presented in this paper. Currently, we are working in that direction.

\subsection{Centrality measures}Although we have identified some pros and cons of using TDA in this particular experiment, a reasonable next step would be to compare our results with those obtained using different methods.  For that purpose, in line with the recent literature on centrality measures presented in the Introduction \cite{huynh2024understanding,nikolaou2020identification,song2017analysis,sun2023centrality}, we have computed four centrality metrics in order to also obtain information of the Spanish Airport Network: 
\begin{enumerate}
    \item \textit{Betweenness centrality:} It measures how often a vertex appears on the shortest paths between other vertices. Vertices with high betweenness can play a key role in a network by influencing the flow of information between other vertices.
    \item \textit{Degree centrality:} It quantifies the ratio between the number of edges links of certain node by the total number of connections of the network.
    \item \textit{Closeness centrality:} It computes the average shortest path from a given vertex to all other vertices in the network.
    \item \textit{Local Reaching centrality:} It indicates how effectively a node can reach others within its local neighborhood. It evaluates the proportion of nodes that can be accessed from the given node within a certain radius, and ultimately selects the maximum value.
\end{enumerate}
For readers interested in deeper definitions and implications of these metrics, we provide the following references: \cite{local_reaching, wasserman1994social}.

We have calculated these metrics for the whole European Airport Network for every day in the Summer Season 2018 (217 days), building a directed and weighted graph for each day: nodes are the airports, direct edges point out that a flight goes from airport $v_1$ to airport $v_2$, and the weight of each edge is the number of flights in that day between airport $v_1$ and $v_2$. Once we built the 217 graphs, we calculated the centrality metrics of each day and, finally the average of each metric, extracting only the Spanish airports results and sorted them accordingly. The whole table with these results can be found in Appendix \ref{app:centrality}. Additionally, Table \ref{table:nodes_edges} provides a summary of the number of edges and nodes in these networks.

\begin{table}[]
\begin{tabular}{|c|c|c|}
\hline \textbf{Num. vertices} & \textbf{Average Edges} & \textbf{Standard Deviation Edges} \\ \hline
1247            & 2160.80       & 185.91\\    \hline              
\end{tabular}
\caption{Summary of the network used to compute centrality measures: the number of vertices remains constant each day, while the second and third columns provide the average and standard deviation of the number of edges in the directed graph of the Spanish airport network.}
\label{table:nodes_edges}
\end{table}

The first clear advantage of using these centrality methods over TDA is the computation time, as TDA is much slower. However, these two methods approach the problem in completely different ways.  Centrality measures, based on graph theory, provide straightforward descriptive insights without involving complex analysis or drawing inferences, whereas TDA tries to focus on interconectedness and searching patterns under different perspectives we previously presented. Although graphs can be enriched with a lot of information and can be weighted in many ways, in this particular study, the information they encapsulate is simpler than the data uniquely created for the footprints of the airports. For example, the geographical relation discovered for LEAS, LEXJ, LEST and LEVX airports are lost with centrality measures.   Moreover, as we previously mentioned, the point cloud used for the persistence diagrams could be made more complex by including additional information, such as economic data of the flights, their order during the day, or even their capacity.  In Limitations' Section, we explore in more detail these possible upgrades.

Table \ref{table:centrality} presents valuable insights. First, Local reaching and Betweenness columns are practically the same.  Although we compute four different aspects of our network, finally we are only basing our new classification over three of them. Second, in terms of the sorting, we also see the Special Group (LEMD, LEBL and LEPA) in the three first positions of the table but, as well as our TDA experiment has detected, Palma's airport has close Degree and Closeness metrics with airports of other groups. Furthermore, for the other atypical group, the Canary's one, we can see in the table that its variability is almost perfectly captured (except for La Palma's airport (GCLA)). Big airports such as Gran Canaria (GCLP) and Tenerife South (GCTS) are in high positions, the big/medium ones are spared around the whole ranking and the smallest ones (La Gomera's (GCGM) and GCHI (GCHI)), with a few flights per week, are the last two of the table.

Additionally, we observe that some high-traffic airports, such as Bilbao (LEBB) and La Palma (GCLA), are far from similar airports. Also, Zaragoza's airport (LEZG) which not operates flights with passengers, is in position thirteen of the ranking, but its uniqueness is not detected as it was in the footprint study. 

Regarding small airports (Group 3), Girona's airport (LEGE) is placed in top positions reflecting its high number of connections and the improvement in number of flights that suffers during Summer seasons. TDA did also reflect this pattern. The remaining airports are in the low part of the table, also agreeing with the actual network classification.

Overall, the classification provided by the centrality measures is largely similar to the one from AENA, with only a few cases suggesting new areas for study that could enhance the classification with more data-driven insights.


We believe that centrality measures are also very powerful and promising for ATM. In our experiment, we were not able to add as much information to the graph as the point clouds used for TDA, but the computation time is considerably shorter than that required for computing the distance matrix $\mathcal{D}$. Both methods have similarities with the current Spanish Airport Network classification and present nuances that enrich it,  being the one provided by centrality measures closer to the actual. We consider that a combination of TDA techniques and centrality measures could be more powerful than using them alone. We are currently working in this direction for future studies.

To sum up, in this section we presented a proof of concept demonstrating how TDA can enhance existing classifications. As we have shown, experiments like the one we conducted can lead to exploring new ways of classifying airports in a more data-driven manner, optimizing their potential and resources. Moreover, TDA can open new perspectives for addressing ATM problems and challenges due to its numerous advantages and its capability to successfully combine with existing techniques.

\section{Limitations}


This study applied Topological Data Analysis (TDA) techniques to analyze the operational behavior of various airports in Spain using data from a single season in 2018. While the results identified clusters of airports with similar behaviors, several limitations affect the generalizability of the findings and suggest potential improvements for future research.

One of the main limitations of this study is that it relies solely on data from a single operational season. This makes it difficult to extrapolate the results to other time periods, as the observed patterns may be influenced by seasonality and interannual variations in airport operations. To draw more representative conclusions, it would be advisable to extend the analysis to include data from multiple years and seasons.

Another aspect that has not been considered is the potential presence of exceptional events on certain days, such as technical, meteorological, or social incidents, which may have altered the usual behaviour of the airports. These events introduce variability into the data, which could skew the results if not appropriately filtered. Therefore, filtering out atypical days would be beneficial in reducing such bias.


Another factor missing from the analysis is contextual information on operational conditions, such as weather, capacity restrictions, and flight demand. Including this data would allow for a more granular analysis and better distinctions between days with similar conditions, enhancing the precision of the findings.


Furthermore, the study did not explore how airports adapt to varying operational conditions. Classifying airports based on their responses to high-demand days or adverse weather could reveal deeper insights into their adaptability and operational resilience.


The analysis would benefit from expanding the dataset to cover multiple years and refining the filtering process to group days with similar characteristics. This would reduce unexplained heterogeneity and improve the quality of the inferences drawn.

Finally, the computational speed of TDA tools is slower compared to other techniques, such as centrality measures. Combining different methods and perspectives, and improving the parallelization of the code would be beneficial for the manageability and development of the process. Although time performance is always important, since we are focused on classification, it is acceptable to sacrifice some computational time in order to capture more nuances.





In conclusion, extending the analysis to include multiple seasons and years, filtering out atypical events, and incorporating supplementary data on airport conditions, as well as improving the speed of computation of the proposal, would enhance the robustness and depth of the study. Classifying airports based on operational scenarios would also provide a more detailed understanding of their behavior, improving the overall accuracy and reliability of the findings.

\section{Conclusions}

Airports and Air Traffic Management Systems are complex sociotechnical structures whose interconnectedness, interdependencies, and complexity generate a huge amount of flight data.
Despite their large potential, flight data are complex, high-dimensional, and sparse datasets, and they are affected by inconsistencies, errors, high levels of variability, multilevel interactions, dynamic changes and high dimensionality, making them very difficult to analyse and exploit.

While existing research works have made significant contributions to the field of analyzing high-dimensional flight trajectory data, there are still limitations that need to be addressed in future research to improve the accuracy, reliability and scalability of these methods.

To overcome these difficulties, this study addresses the use of Topological Data Analysis (TDA) to analyse flight trajectory data. In particular, as a proof of concept, we try to determine relationships between different variables involved in the spatial and temporal flight trajectory and delays to identify common patterns and anomalies in airport operation, extracting a \textit{footprint} of each airport as a way to identify it. In future studies, it could help to recognize the underlying causes of delays and develop more effective strategies for reducing them and produce, combined to other method and techniques, a more exhaustive classification of groups of airports.

Real-world data from the Spanish network of airports in the Summer Season of 2018 were used for the assessment. For each set of airports, the average persistence landscape has been computed for each airport in that group, based on which we construct a distance matrix $\mathcal{D}$, where each element $d_{ij}$ is the supremum distance between the persistence average landscapes of airport $i$ and airport $j$. After this, we performed a block permutation test in order to check the statistical significance of our experiment. Finally, we studied either $\mathcal{D}$ and the $p$-values in order to extract some conclusions and compare them with the current Spanish Airport Network classification and another non-TDA method based on centrality measures. 

The accomplished analysis points out several conclusions regarding the grouping of airports attending to operational criteria, and regarding the identification of patterns and anomalies:

\begin{itemize}
    \item The results of the analysis show how different airport groups follow a sort of cluster, even though there are few in number in each group. The analysis show that the imposed upper limit on Group 3 compelled the incorporation of airports of diverse types. Consequently, Group 3 has been subdivided into four specific subgroups, namely general aviation airports, air bases accessible to civilian traffic, airports with low traffic, and helipads. Likewise, the Canary Group presents a notable variation in the sizes of its airports, ranging from La Gomera with two to four flights per day to Gran Canaria with over one hundred flights per day.
    \item It also allows the identification of airports that clearly differ from their preassigned group (such as Zaragoza’s airport, which differs from the rest of the airports of Group 2). It highlight the Zaragoza’s airport uniqueness, as the only airport in Group 2 that does not facilitate passenger aviation. 
    \item It also helps detecting when an airport is isolated and far from any other airport in the network, such as Adolfo Suarez Madrid-Barajas and Josep Tarradellas Barcelona-El Prat’s airports are from any other airport in Spain. The Average Persistence landscape of these airports also exhibit a big the difference in magnitude of the axes.
    \item Moreover, the diversity in a category of airports is also reflected in the Isomap representation, as in the case of the Canary Airports Group. It is very interesting to check using distance matrix information how the Canary Group behaves with the rest of Spanish airport net.  It seems that its geographical location and the length of its flights (on longer flights the deviation between trajectories increases) are an important 	factor to characterize this group.
    \item It illustrated also that point cloud create based on deviation of trajectories and delays, encodes information regarding geographical location of the airports and origins and destination of its flights. This explain that airports located in close proximity tend to have similar geometric properties, as it is the case with Asturias (LEAS), Santander (LEXJ), Santiago (LEST), and Vigo (LEVX). TDA analysis is capable of capturing that as it can check in the distances of airports with those similar characteristics.
    \item Additionally, it presents some relevant questions as the overlapping could be interpreting as whether some airports can belong to different groups. In particular, we can see that the point corresponding to Zaragoza’s airport  is far from the rest of the airports in Group 2. These outcomes shows how TDA may guide us toward exploring novel approaches to categorize airports in a more data-driven manner, aiming to optimize their full potential and resources. 
    \item Lastly, the $p$-values shown in Appendix \ref{app:p-value} reinforce some of our conclusions and highlight weaknesses in our analysis, indicating future directions for improving our classification.
\end{itemize}

This study has become a proof of concept of how TDA can become a powerful analytical technique to help overcome some of the limitations of existing research in the field of analyzing high-dimensional flight trajectory data for the identification of common traffic patterns in airports.
In particular, during the study, the following advantages of TDA were proposed: 

\begin{itemize}
    \item It can assist in identifying underlying structures and patterns within the data in a manner that is more robust to noise and outliers.
    \item It can also help address the issue of limited data availability by allowing for the integration of different data sources and the extraction of meaningful insights from incomplete or noisy data.
    This can be a crucial tool for solving the problem of scalability by allowing the analysis of large and complex datasets using parallel computing techniques. 
    \item It can attempt to address the lack of standardization by offering a flexible and adaptable framework that is applicable to a diverse array of data formats and types.
\end{itemize}

To the best of our knowledge, apart from some purely introductory contributions \cite{maxdetection,maxtda}, this paper presents the first effective application of TDA to the aerospace field. In particular, no rigorous work had been performed by applying this method to a large amount of aircraft trajectory data in an attempt to anticipate and identify anomalies in aircraft space/time trajectories, infer patterns of behavior at different airports, and classify and characterize airports depending on the distribution of their daily flights via trajectory deviation and delay. This pioneering effort not only expands the scope of TDA but also lays the groundwork for future research in leveraging this method to enhance our understanding of aircraft dynamics, ultimately contributing to improved aviation safety and efficiency.

\clearpage

\appendix

\section{Tables}{\label{tablas}}

\begin{table}[h]
    \centering
\resizebox{350pt}{!}{%
\begin{tabular}{l|cccccc}
{} & LEAB & LEBZ & LELN & LESA & LESB & LEVD \\
\midrule
LEAB & 0.00 & 106.28 & 58.60 & 129.64 & 112.24 & 107.70 \\
LEBZ & 106.28 & 0.00 & 53.22 & 54.47 & 23.73 & 42.73 \\
LELN & 58.60 & 53.22 & 0.00 & 100.65 & 69.99 & 49.32 \\
LESA & 129.64 & 54.47 & 100.65 & 0.00 & 38.71 & 74.42 \\
LESB & 112.24 & 23.73 & 69.99 & 38.71 & 0.00 & 61.62 \\
LEVD & 107.70 & 42.73 & 49.32 & 74.42 & 61.62 & 0.00 \\
\end{tabular}}
    \caption{Distance matrix corresponding to Group 3 - \textit{Air bases open to civilian traffic}}
    \label{tablagrupo3ba}
\end{table}

\begin{table}[h]
    \centering
\resizebox{350pt}{!}{%
\begin{tabular}{l|cccccc}
{} & LEAB & LEBZ & LELN & LESA & LESB & LEVD \\
\midrule
LECU & 123.04 & 34.12 & 85.25 & 27.17 & 16.24 & 59.85 \\
LELL & 162.93 & 75.88 & 125.81 & 35.15 & 57.94 & 94.85 \\
\end{tabular}}
    \caption{Distance matrix corresponding to Group 3 - \textit{Air bases open to civilian traffic} and the Group 3 - \textit{General aviation airports}}
    \label{tablagrupo3ba3g}
\end{table}

\begin{table}[h]
    \centering
\resizebox{350pt}{!}{%
\begin{tabular}{l|cccccc}
{} & LEAB & LEBZ & LELN & LESA & LESB & LEVD \\
\midrule
GEML & 181.95 & 118.74 & 160.00 & 64.59 & 102.96 & 138.88 \\
LEBA & 119.56 & 32.91 & 81.47 & 33.18 & 14.24 & 72.80 \\
LEBG & 144.18 & 58.86 & 109.44 & 16.39 & 40.72 & 83.15 \\
LEGE & 77.39 & 46.06 & 46.60 & 65.33 & 59.18 & 36.74 \\
LEHC & 126.42 & 83.75 & 115.76 & 32.04 & 69.63 & 103.73 \\
LEPP & 102.69 & 17.30 & 56.70 & 43.97 & 34.28 & 34.95 \\
LERJ & 126.29 & 58.37 & 102.48 & 15.52 & 41.34 & 78.15 \\
LESO & 132.98 & 56.64 & 102.03 & 30.05 & 39.73 & 76.50 \\
LEVT & 99.37 & 38.94 & 60.45 & 60.09 & 53.85 & 32.98 \\
\end{tabular}}
    \caption{Distance matrix corresponding to Group 3 - \textit{Air bases open to civilian traffic} and the Group 3 - \textit{Airports with low traffic}}
    \label{tablagrupo3ba3bt}
\end{table}

\begin{table}[h]
    \centering
\resizebox{350pt}{!}{%
\begin{tabular}{l|cccccc}
{} & LEAB & LEBZ & LELN & LESA & LESB & LEVD \\
\midrule
LEAM & 112.78 & 48.08 & 72.95 & 67.41 & 62.13 & 35.08 \\
LEAS & 118.22 & 39.16 & 75.56 & 57.19 & 53.91 & 42.99 \\
LECO & 139.56 & 39.85 & 93.06 & 44.68 & 42.57 & 54.90 \\
LEGR & 131.81 & 39.66 & 88.52 & 55.22 & 52.58 & 49.40 \\
LEJR & 98.17 & 42.70 & 52.22 & 61.72 & 57.50 & 26.37 \\
LERS & 109.29 & 48.25 & 58.51 & 68.06 & 62.81 & 25.06 \\
LEVX & 120.37 & 35.99 & 77.80 & 50.96 & 49.25 & 41.39 \\
LEXJ & 118.34 & 41.27 & 69.07 & 56.62 & 55.20 & 30.72 \\
LEZG & 125.63 & 189.13 & 143.81 & 192.98 & 189.22 & 191.02 \\
\end{tabular}}
    \caption{Distance matrix corresponding to Group 3 - \textit{Air bases open to civilian traffic} and the Group 2}
    \label{tablagrupo3ba2}
\end{table}

\begin{table}[h]
    \centering
\resizebox{350pt}{!}{%
\begin{tabular}{l|cccccc}
{} & LEAB & LEBZ & LELN & LESA & LESB & LEVD \\
\midrule
LEAL & 47.75 & 66.99 & 35.10 & 82.37 & 67.93 & 69.30 \\
LEBB & 84.02 & 36.38 & 50.21 & 50.66 & 46.71 & 37.36 \\
LEIB & 84.27 & 42.43 & 47.18 & 58.95 & 46.47 & 44.49 \\
LEMG & 41.04 & 104.21 & 57.36 & 119.10 & 105.47 & 106.33 \\
LEMH & 119.76 & 37.16 & 69.92 & 53.32 & 50.55 & 31.70 \\
LEST & 102.49 & 42.78 & 59.97 & 56.79 & 54.88 & 27.46 \\
LEVC & 55.67 & 54.35 & 34.99 & 76.36 & 56.65 & 56.30 \\
LEZL & 87.26 & 36.69 & 46.76 & 55.46 & 47.77 & 29.80 \\
\end{tabular}}
    \caption{Distance matrix corresponding to Group 3 - \textit{Air bases open to civilian traffic} and the Group 1}
    \label{tablagrupo3ba1}
\end{table}

\begin{table}[h]
    \centering
\resizebox{350pt}{!}{%
\begin{tabular}{l|cccccc}
{} & LEAB & LEBZ & LELN & LESA & LESB & LEVD \\
\midrule
GCFV & 46.31 & 78.20 & 46.51 & 98.36 & 81.18 & 80.42 \\
GCGM & 182.02 & 122.09 & 161.68 & 68.55 & 107.16 & 142.10 \\
GCHI & 184.57 & 124.98 & 164.57 & 71.30 & 109.97 & 145.08 \\
GCLA & 51.30 & 57.61 & 38.93 & 100.42 & 71.33 & 56.45 \\
GCLP & 45.91 & 79.78 & 48.25 & 86.12 & 82.01 & 79.44 \\
GCRR & 50.11 & 63.25 & 35.79 & 82.53 & 63.65 & 65.42 \\
GCTS & 88.64 & 143.24 & 101.00 & 142.02 & 144.74 & 143.84 \\
GCXO & 112.29 & 34.75 & 73.34 & 38.31 & 39.44 & 44.84 \\
\end{tabular}}
    \caption{Distance matrix corresponding to Group 3 - \textit{Air bases open to civilian traffic} and the Canary Group}
    \label{tablagrupo3bac}
\end{table}

\begin{table}[h]
    \centering
\resizebox{350pt}{!}{%
\begin{tabular}{l|cccccc}
{} & LEAB & LEBZ & LELN & LESA & LESB & LEVD \\
\midrule
LEBL & 306.99 & 353.02 & 315.33 & 343.49 & 356.43 & 350.35 \\
LEMD & 453.15 & 462.99 & 457.91 & 460.78 & 465.56 & 462.91 \\
LEPA & 72.94 & 33.36 & 30.60 & 70.09 & 44.31 & 35.09 \\
\end{tabular}}
    \caption{Distance matrix corresponding to Group 3 - \textit{Air bases open to civilian traffic} and the Special Group}
    \label{tablagrupo3bae}
\end{table}

\begin{table}[h]
    \centering
\resizebox{150pt}{!}{\begin{tabular}{l|cc}
{} &       LECU &       LELL \\
\midrule
LECU &   0.00 &  41.91 \\
LELL &  41.91 &   0.00 \\
\end{tabular}}
    \caption{Distance matrix corresponding to Group 3 - \textit{General aviation airports}}
    \label{tablagrupo3g}
\end{table}

\begin{table}[h]
    \centering
\resizebox{\columnwidth}{!}{%
\begin{tabular}{l|ccccccccc}
{} & GEML & LEBA & LEBG & LEGE & LEHC & LEPP & LERJ & LESO & LEVT \\
\midrule
LECU & 89.06 & 21.01 & 28.24 & 54.45 & 59.20 & 30.10 & 26.73 & 23.76 & 48.71 \\
LELL & 47.78 & 47.37 & 21.46 & 85.80 & 64.00 & 69.90 & 50.25 & 32.80 & 67.38 \\
\end{tabular}}
    \caption{Distance matrix corresponding to Group 3 - \textit{General aviation airports} and Group 3 - \textit{Airports with low traffic}}
    \label{tablagrupo3g3bt}
\end{table}

\begin{table}[h]
    \centering
\resizebox{\columnwidth}{!}{%
\begin{tabular}{l|ccccccccc}
{} & LEAM & LEAS & LECO & LEGR & LEJR & LERS & LEVX & LEXJ & LEZG \\
\midrule
LECU & 57.34 & 48.19 & 34.79 & 46.54 & 52.02 & 58.10 & 42.28 & 49.27 & 191.09 \\
LELL & 74.92 & 63.95 & 49.94 & 61.65 & 75.69 & 75.95 & 56.54 & 66.50 & 227.78 \\
\end{tabular}}
    \caption{Distance matrix corresponding to Group 3 - \textit{General aviation airports} and Group 2}
    \label{tablagrupo3g2}
\end{table}

\begin{table}[h]
    \centering
\resizebox{450pt}{!}{%
\begin{tabular}{l|cccccccc}
{} & LEAL & LEBB & LEIB & LEMG & LEMH & LEST & LEVC & LEZL \\
\midrule
LECU & 76.43 & 40.00 & 44.15 & 113.35 & 44.85 & 48.73 & 68.11 & 41.42 \\
LELL & 116.14 & 79.10 & 81.42 & 152.04 & 65.81 & 67.53 & 108.09 & 79.07 \\
\end{tabular}}
    \caption{Distance matrix corresponding to Group 3 - \textit{General aviation airports} and Group 1}
    \label{tablagrupo3g1}
\end{table}

\begin{table}[h]
    \centering
\resizebox{450pt}{!}{%
\begin{tabular}{l|cccccccc}
{} & GCFV & GCGM & GCHI & GCLA & GCLP & GCRR & GCTS & GCXO \\
\midrule
LECU & 91.67 & 94.84 & 97.42 & 85.22 & 80.52 & 73.31 & 139.60 & 27.44 \\
LELL & 131.56 & 59.15 & 60.76 & 124.78 & 119.50 & 113.23 & 176.94 & 52.46 \\
\end{tabular}}
    \caption{Distance matrix corresponding to Group 3 - \textit{General aviation airports} and the Canary Group}
    \label{tablagrupo3gC}
\end{table}

\begin{table}[h]
    \centering
\resizebox{200pt}{!}{\begin{tabular}{l|ccc}
{} &        LEBL &        LEMD &       LEPA \\
\midrule
LECU &  344.80 &  459.58 &  56.15 \\
LELL &  367.96 &  465.56 &  96.01 \\
\end{tabular}}
    \caption{Distance matrix corresponding to Group 3 - \textit{General aviation airports} and the Special Group}
    \label{tablagrupo3ge}
\end{table}

\begin{table}[h]
    \centering
\resizebox{\columnwidth}{!}{%
\begin{tabular}{l|ccccccccc}
{} & LEAM & LEAS & LECO & LEGR & LEJR & LERS & LEVX & LEXJ & LEZG \\
\midrule
GCFV & 81.48 & 88.23 & 108.45 & 100.62 & 67.72 & 80.53 & 89.26 & 87.79 & 117.50 \\
GCGM & 122.74 & 109.65 & 101.75 & 110.15 & 124.66 & 121.86 & 106.46 & 119.41 & 242.30 \\
GCHI & 125.27 & 111.88 & 104.02 & 112.59 & 127.40 & 124.45 & 109.13 & 122.06 & 244.70 \\
GCLA & 74.32 & 75.58 & 96.38 & 90.27 & 53.69 & 61.72 & 79.79 & 73.31 & 149.05 \\
GCLP & 72.40 & 85.90 & 100.04 & 92.57 & 70.58 & 78.18 & 79.96 & 80.55 & 112.87 \\
GCRR & 63.34 & 70.80 & 90.06 & 82.17 & 49.86 & 63.53 & 70.86 & 69.64 & 127.64 \\
GCTS & 129.36 & 149.25 & 161.29 & 153.26 & 131.04 & 142.48 & 140.05 & 143.84 & 61.07 \\
GCXO & 41.27 & 40.97 & 45.53 & 34.54 & 33.13 & 42.15 & 30.55 & 34.00 & 177.50 \\
\end{tabular}}
    \caption{Distance matrix corresponding to Group 2 and the Canary Group}
    \label{tablagrupo2C}
\end{table}

\begin{table}[h]
    \centering
\resizebox{\columnwidth}{!}{%
\begin{tabular}{l|ccccccccc}
{} & LEAM & LEAS & LECO & LEGR & LEJR & LERS & LEVX & LEXJ & LEZG \\
\midrule
LEBL & 345.30 & 360.04 & 365.06 & 354.35 & 343.12 & 349.85 & 348.75 & 351.07 & 185.01 \\
LEMD & 465.38 & 465.56 & 465.56 & 461.14 & 464.49 & 464.01 & 463.71 & 465.56 & 378.29 \\
LEPA & 45.21 & 46.44 & 68.64 & 61.62 & 26.45 & 37.20 & 50.77 & 46.78 & 160.32 \\
\end{tabular}}
    \caption{Distance matrix corresponding to Group 2 and the Special Group}
    \label{tablagrupo2E}
\end{table}

\begin{table}[h]
    \centering
\resizebox{450pt}{!}{%
\begin{tabular}{l|cccccccc}
{} & LEAL & LEBB & LEIB & LEMG & LEMH & LEST & LEVC & LEZL \\
\midrule
GCFV & 26.65 & 52.65 & 52.89 & 26.04 & 90.16 & 71.77 & 27.20 & 55.99 \\
GCGM & 139.67 & 117.15 & 126.72 & 174.18 & 117.29 & 118.14 & 131.50 & 120.68 \\
GCHI & 142.64 & 119.99 & 129.74 & 176.85 & 120.03 & 120.86 & 134.19 & 123.73 \\
GCLA & 23.48 & 50.21 & 46.73 & 55.59 & 73.78 & 60.72 & 33.89 & 46.31 \\
GCLP & 43.73 & 57.89 & 48.47 & 33.18 & 87.13 & 65.61 & 46.17 & 60.46 \\
GCRR & 12.96 & 34.19 & 34.61 & 42.75 & 72.65 & 53.39 & 15.37 & 37.72 \\
GCTS & 89.16 & 108.25 & 102.16 & 62.74 & 150.96 & 129.57 & 97.67 & 115.63 \\
GCXO & 65.62 & 28.50 & 31.22 & 101.74 & 41.14 & 27.03 & 57.35 & 30.71 \\
\end{tabular}}
    \caption{Distance matrix corresponding to Group 1 and the Canary Group}
    \label{tablagrupo1C}
\end{table}

\begin{table}[h]
    \centering
\resizebox{450pt}{!}{%
\begin{tabular}{l|cccccccc}
{} & GCFV & GCGM & GCHI & GCLA & GCLP & GCRR & GCTS & GCXO \\
\midrule
LEBL & 296.47 & 366.86 & 367.96 & 327.56 & 274.82 & 306.85 & 219.22 & 319.87 \\
LEMD & 448.45 & 465.56 & 465.56 & 463.43 & 419.33 & 460.48 & 436.09 & 449.87 \\
LEPA & 46.63 & 132.85 & 135.83 & 30.78 & 67.96 & 34.24 & 119.95 & 43.95 \\
\end{tabular}}
    \caption{Distance matrix corresponding to the Canary Group and the Special Group}
    \label{tablagrupoCE}
\end{table}

\clearpage
\afterpage{
\begin{landscape}
\pagestyle{empty}%
\section{Statistical significance: p-values}{\label{app:p-value}}
\topskip0pt
\begin{table}[htb!]
\resizebox{500pt}{!}{
\begin{tabular}{l||ll:llllll:lllllllll}
     & LECU   & LELL   & LEAB   & LEBZ   & LELN   & LESA   & LESB   & LEVD   & GEML   & LEBG   & LEBA   & LEGE   & LEHC   & LERJ   & LEPP   & LESO   & LEVT   \\ \hline \hline
LECU &        &        &        &    &        & 0.0675 & 0.5575 &        &        & 0.1025 & 0.68   &        &    & 0.13   &    & 0.08   &        \\ 
LELL &        &        &        &        &        &        &        &        &        & 0.085  &        &        &        &        &        &        &        \\ \hdashline
LEAB &        &        &        &        &  &        &        &        &        &        &        &        &        &        &        &        &        \\
LEBZ &    &        &        &        &        &        & 0.1875 &        &        &        & 0.255  &        &        &        & 0.3425 &        &    \\
LELN &        &        &  &        &        &        &        &        &        &        &        &        &        &        &        &        &        \\
LESA & 0.0675 &        &        &        &        &        &  &        &        & 0.56   & 0.2225 &        & 0.3575 & 0.6925 &        &  &        \\
LESB & 0.5575 &        &        & 0.1875 &        &  &        &        &        &    & 0.8825 &        &  &   &   &        &        \\
LEVD &        &        &        &        &        &        &        &        &        &        &        &    &    &        &   &        &  \\ \hdashline
GEML &        &        &        &        &        &        &        &        &        &        &        &        &  &        &        &        &        \\
LEBG & 0.1025 & 0.085  &        &        &        & 0.56   &    &        &        &        & 0.1975 &        & 0.1    & 0.155  &        &        &        \\
LEBA & 0.68   &        &        & 0.255  &        & 0.2225 & 0.8825 &        &        & 0.1975 &        &        &   & 0.29   & 0.06   & 0.1725 &    \\
LEGE &        &        &        &        &        &        &        &    &        &        &        & 0.1225 &        &   &        & 0.25   \\
LEHC &    &        &        &        &        & 0.3575 &  &    &  & 0.1    &   & 0.1225 &        & 0.37   & 0.18   & 0.18   &        \\
LERJ & 0.13   &        &        &        &        & 0.6925 &   &        &        & 0.155  & 0.29   &        & 0.37   &        &        &  &        \\
LEPP &    &        &        & 0.3425 &        &        &   &   &        &        & 0.06   &   & 0.18   &        &        &        & 0.155  \\
LESO & 0.08   &        &        &        &        &  &        &        &        &        & 0.1725 &        & 0.18   &        &        &        &        \\
LEVT &        &        &        &    &        &        &        &  &        &        &        & 0.25   &        &        & 0.155  &        &        \\ \hline
LECO &        &        &        &        &        &        &        &        &        &        &        &        &   &        &        &  &  \\
LEAM &        &        &        &        &        &        &        &        &        &        &        &    &   &        &    &        & 0.505  \\
LEAS &        &        &        &        &        &        &        &        &        &        &        &        &   &        &   &        & 0.18   \\
LEGR &        &        &        &        &        &        &        &        &        &        &        &        &  &        &        &        &   \\
LEJR &        &        &        &        &        &        &        & 0.08   &        &        &        & 0.2225 &    &        &    &        & 0.6125 \\
LERS &        &        &        &        &        &        &        &    &        &        &        &    &  &        &        &        & 0.5025 \\
LEXJ &        &        &        &        &        &        &        &  &        &        &        &        & 0.145  &        &  &        & 0.3125 \\
LEST &        &        &        &        &        &        &        & 0.065  &        &        &        & 0.1475 &   &        & 0.0875 &        & 0.8    \\
LEVX &        &        &        &        &        &        &        &        &        &        &        &        &   &        & 0.13   &        & 0.2575 \\
LEZG &        &        &        &        &        &        &        &        &        &        &        &        &  &        &        &        &        \\ \hline
LEAL &        &        & 0.185  &        & 0.1025 &        &        &        &        &        &        & 0.15   &    &        &        &        &        \\
LEBB &  &        &        &  &    &        &        &  &        &        & 0.08   & 0.25   & 0.2475 &        & 0.065  &        & 0.38   \\
LEIB &        &        &        &    &   &        &  &        &        &        & 0.2575 & 0.415  & 0.245  &        &  &        & 0.1925 \\
LEMG &        &        & 0.3525 &        &    &        &        &        &        &        &        &        &  &        &        &        &        \\
LEMH &        &        &        &        &        &        &        &        &        &        &        &        & &        & 0.0975 &        & 0.16   \\
LEZL &        &        &        &        &        &        &    & 0.135  &        &        & 0.065  & 0.6775 & 0.2275 &        & 0.2275 &        & 0.5825 \\
LEVC &        &        & 0.0825 &        & 0.1025 &        &        &        &        &        &  & 0.375  &        &        &        &        &  \\ \hline
GCFV &        &        & 0.2225 &        &    &        &        &        &        &        &        &   &  &        &        &        &        \\
GCGM &        &        &        &        &        &        &        &        &        &        &        &        &        &        &        &        &        \\
GCHI &        &        &        &        &        &        &        &        &        &        &        &        &   &        &        &        &        \\
GCLA &        &        & 0.0775 &        &   &        &        &        &        &        &        &        &        &        &        &        &        \\
GCLP &        &        & 0.325  &        &  &        &        &        &        &        &        &        &  &        &        &        &        \\
GCRR &        &        & 0.1025 &        & 0.1075 &        &        &        &        &        &        & 0.1175 &  &        &        &        &   \\
GCTS &        &        &  &        &        &        &        &        &        &        &        &        & 0.33   &        &        &        &        \\
GCXO & 0.1125 &        &        &    &        &    &  &        &        &    & 0.21   & 0.0525 & 0.1575 &    & 0.13   & 0.0625 & 0.13   \\ \hline
LEMD &        &        &        &        &        &        &        &        &        &        &        &        & 0.2375 &        &        &        &        \\
LEBL &        &        &        &        &        &        &        &        &        &        &        &        &  &        &        &        &        \\
LEPA &        &        &        &   & 0.14   &        &        &        &        &        & 0.0875 & 0.2725 & 0.12   &        &   &        & 0.07  
\end{tabular}}
\caption{$p$-values between Group 3 and the rest of the airports}
\label{table:p-valuegroup3}
\end{table}

\begin{table}[htb!]
\resizebox{520pt}{!}{
\begin{tabular}{l||llllllllll|lllllll}
     & LECO   & LEAM   & LEAS   & LEGR   & LEJR   & LERS   & LEXJ   & LEST   & LEVX   & LEZG   & LEAL   & LEBB   & LEIB   & LEMG   & LEMH   & LEZL   & LEVC   \\ \hline \hline
LECU &        &        &        &        &        &        &        &        &        &        &        &  &        &        &        &        &        \\
LELL &        &        &        &        &        &        &        &        &        &        &        &        &        &        &        &        &        \\ \hdashline
LEAB &        &        &        &        &        &        &        &        &        &        & 0.185  &        &        & 0.3525 &        &        & 0.0825 \\
LEBZ &        &        &        &        &        &        &        &        &        &        &        &  &        &        &        &        &        \\
LELN &        &        &        &        &        &        &        &        &        &        & 0.1025 &        &        &        &        &        & 0.1025 \\
LESA &        &        &        &        &        &        &        &        &        &        &        &        &        &        &        &        &        \\
LESB &        &        &        &        &        &        &        &        &        &        &        &        &  &        &        &        &        \\
LEVD &        &  &        &        & 0.08   &        &  & 0.065  &        &        &        &  &        &        &        & 0.135  &        \\ \hdashline
GEML &        &        &        &        &        &        &        &        &        &        &        &        &        &        &        &        &        \\
LEBG &        &        &        &        &        &        &        &        &        &        &        &        &  &        &        &        &        \\
LEBA &        &        &        &        &        &        &        &        &        &        &        & 0.08   & 0.2575 &        &        & 0.065  &  \\
LEGE &        &   &        &        & 0.2225 &        &        & 0.1475 &        &        & 0.15   & 0.25   & 0.415  &        &        & 0.6775 & 0.375  \\
LEHC &   &   &   &  &    &  & 0.145  &   &   &  &    & 0.2475 & 0.245  &  &  & 0.2275 &  \\
LERJ &        &        &        &        &        &        &        &        &        &        &        &        &        &        &        &        &        \\
LEPP &        &        &        &        &        &        &  & 0.0875 & 0.13   &        &        & 0.065  &  &        & 0.0975 & 0.2275 &        \\
LESO &        &        &        &        &        &        &        &        &        &        &        &        &        &        &        &        &        \\
LEVT &        & 0.505  & 0.18   &   & 0.6125 & 0.5025 & 0.3125 & 0.8    & 0.2575 &        &        & 0.38   & 0.1925 &        & 0.16   & 0.5825 &  \\ \hline
LECO &        &   & 0.075  & 0.2725 &        &        & 0.0775 &        & 0.1325 &        &        &        &        &        &        &        &        \\
LEAM &   &        & 0.1875 & 0.055  & 0.275  & 0.3925 & 0.425  & 0.3125 &  &        &        & 0.1    & 0.0875 &        & 0.1325 & 0.14   &   \\
LEAS & 0.075  & 0.1875 &        & 0.31   & 0.105  & 0.18   & 0.4075 & 0.1225 & 0.575  &        &        &        &        &        & 0.42   &  &        \\
LEGR & 0.2725 & 0.055  & 0.31   &        &        &        & 0.1575 &  & 0.52   &        &        &        &        &        & 0.1575 &        &        \\
LEJR &        & 0.275  & 0.105  &        &        & 0.5175 & 0.2275 & 0.4625 & 0.0575 &        &        & 0.2725 & 0.1325 &        & 0.14   & 0.6275 &    \\
LERS &        & 0.3925 & 0.18   &        & 0.5175 &        & 0.3075 & 0.2325 &        &        &        &  &        &        & 0.225  & 0.1375 &        \\
LEXJ & 0.0775 & 0.425  & 0.4075 & 0.1575 & 0.2275 & 0.3075 &        & 0.415  & 0.3625 &        &        &        &        &        & 0.7725 & 0.08   &        \\
LEST &        & 0.3125 & 0.1225 &  & 0.4625 & 0.2325 & 0.415  &        & 0.285  &        &        & 0.2425 & 0.12   &        & 0.115  & 0.2775 &    \\
LEVX & 0.1325 &  & 0.575  & 0.52   & 0.0575 &        & 0.3625 & 0.285  &        &        &        &  &        &        & 0.17   & 0.0725 &        \\
LEZG &        &        &        &        &        &        &        &        &        &        &        &        &        &        &        &        &        \\ \hline
LEAL &        &        &        &        &        &        &        &        &        &        &        & 0.085  & 0.13   & 0.1125 &        &    & 0.6825 \\
LEBB &        & 0.1    &        &        & 0.2725 &        &        & 0.2425 &  &        & 0.085  &        & 0.61   &        &        & 0.2675 & 0.19   \\
LEIB &        & 0.0875 &        &        & 0.1325 &        &        & 0.12   &        &        & 0.13   & 0.61   &        &        &        & 0.5175 & 0.2075 \\
LEMG &        &        &        &        &        &        &        &        &        &        & 0.1125 &        &        &        &        &        &        \\
LEMH &    & 0.1325 & 0.42   & 0.1575 & 0.14   & 0.225  & 0.7725 & 0.115  & 0.17   &        &        &        &        &        &    &        \\
LEZL &        & 0.14   &  &        & 0.6275 & 0.1375 & 0.08   & 0.2775 & 0.0725 &        &    & 0.2675 & 0.5175 &        &    &        & 0.145  \\
LEVC &        &   &        &        &        &        &        &    &        &        & 0.6825 & 0.19   & 0.2075 &        &        & 0.145  &        \\ \hline
GCFV &        &        &        &        &        &        &        &        &        &        & 0.3325 &    &    & 0.415  &        &        & 0.2875 \\
GCGM &        &        &        &        &        &        &        &        &        &        &        &        &        &        &        &        &        \\
GCHI &        &        &        &        &        &        &        &        &        &        &        &        &        &        &        &        &        \\
GCLA &        &        &        &        &        &        &        &        &        &        & 0.3625 &  &  &    &        &  & 0.14   \\
GCLP &        &        &        &        &        &        &        &        &        &        & 0.0625 &        & 0.0575 & 0.2675 &        &        &    \\
GCRR &        &        &        &        &   &        &        &        &        &        & 0.8525 & 0.12   & 0.1    &  &        & 0.0825 & 0.7325 \\
GCTS &        &        &        &        &        &        &        &        &        & 0.095  &        &        &        &  &        &        &        \\
GCXO &  &   &  &    & 0.0875 &   &  & 0.2375 & 0.1    &        &        & 0.245  & 0.235  &        &  & 0.1875 &        \\ \hline
LEMD &        &        &        &        &        &        &        &        &        &        &        &        &        &        &        &        &        \\
LEBL &        &        &        &        &        &        &        &        &        &        &        &        &        &        &        &        &        \\
LEPA &        &        &        &        & 0.1075 &        &        &        &        &        & 0.0675 & 0.3575 & 0.3325 &        &        & 0.6025 & 0.3075
\end{tabular}}
\caption{$p$-values between Group 2 and Group 1 and the rest of the airports}
\label{table:p-valuegroups12}

\end{table}

\begin{table}[htb!]
\resizebox{400pt}{!}{
\begin{tabular}{l||llllllll|lll}
     & GCFV   & GCGM & GCHI & GCLA   & GCLP   & GCRR   & GCTS   & GCXO   & LEMD   & LEBL   & LEPA   \\ \hline \hline
LECU &        &      &      &        &        &        &        & 0.1125 &        &        &        \\
LELL &        &      &      &        &        &        &        &        &        &        &        \\ \hdashline
LEAB & 0.2225 &      &      & 0.0775 & 0.325  & 0.1025 &  &        &        &        &        \\
LEBZ &        &      &      &        &        &        &        &    &        &        &        \\
LELN &    &      &      &   &  & 0.1075 &        &        &        &        & 0.14   \\
LESA &        &      &      &        &        &        &        &    &        &        &        \\
LESB &        &      &      &        &        &        &        &  &        &        &        \\
LEVD &        &      &      &        &        &        &        &  &        &        &  \\ \hdashline
GEML &        &      &      &        &        &        &        &        &        &        &        \\
LEBG &        &      &      &        &        &        &        &    &        &        &        \\
LEBA &        &      &      &        &  &        &        & 0.21   &        &        & 0.0875 \\
LEGE &        &      &      &        &        & 0.1175 &        & 0.0525 &        &        & 0.2725 \\
LEHC &  &      &  &        &  &  & 0.33   & 0.1575 & 0.2375 &        &    \\
LERJ &        &      &      &        &        &        &        &    &        &        &        \\
LEPP &        &      &      &        &        &        &        & 0.13   &        &        & 0.075  \\
LESO &        &      &      &        &        &        &        & 0.0625 &        &        &        \\
LEVT &        &      &      &        &        &        &        & 0.13   &        &        & 0.07   \\ \hline
LECO &        &      &      &        &        &        &        &        &        &        &        \\
LEAM &        &      &      &        &        &        &        &   &        &        &        \\
LEAS &        &      &      &        &        &        &        &  &        &        &        \\
LEGR &        &      &      &        &        &        &        &    &        &        &        \\
LEJR &        &      &      &        &        &        &        & 0.0875 &        &        & 0.1075 \\
LERS &        &      &      &        &        &        &        &   &        &        &  \\
LEXJ &        &      &      &        &        &        &        &  &        &        &  \\
LEST &        &      &      &        &        &        &        & 0.2375 &        &        & 0.06   \\
LEVX &        &      &      &        &        &        &        & 0.1    &        &        &        \\
LEZG &        &      &      &        &        &        & 0.095  &        &        &        &        \\ \hline
LEAL & 0.3325 &      &      & 0.3625 & 0.0625 & 0.8525 &        &        &        &        & 0.0675 \\
LEBB &    &      &      &  &        & 0.12   &        & 0.245  &        &        & 0.3575 \\
LEIB &    &      &      &        & 0.0575 & 0.1    &        & 0.235  &        &        & 0.3325 \\
LEMG & 0.415  &      &      &        & 0.2675 &  &  &        &        &        &        \\
LEMH &        &      &      &        &        &        &        &  &        &        &        \\
LEZL &        &      &      &  &        &        &        & 0.1875 &        &        & 0.6025 \\
LEVC & 0.2875 &      &      & 0.14   &    & 0.7325 &        &        &        &        & 0.3075 \\ \hline
GCFV &        &      &      & 0.215  & 0.145  & 0.3025 &  &        &        &        &    \\
GCGM &        &      &      &        &        &        &        &        &        &        &        \\
GCHI &        &      &      &        &        &        &        &        &        &        &        \\
GCLA & 0.215  &      &      &        &        & 0.3925 &        &        &        &        & 0.14   \\
GCLP & 0.145  &      &      &        &        &        &    &        &        &        &        \\
GCRR & 0.3025 &      &      & 0.3925 &        &        &        &   &        &        & 0.0825 \\
GCTS &  &      &      &        &    &        &        &        &        &        &        \\
GCXO &        &      &      &        &        &        &        &        &        &        &  \\ \hline
LEMD &        &      &      &        &        &        &        &        &  &        &        \\
LEBL &        &      &      &        &        &        &        &        &  &        &        \\
LEPA &        &      &      & 0.14   &        & 0.0825 &        &  &        &        &       
\end{tabular}}
\caption{p-values between Canary Group and Special Group and the rest of the airports}
\label{table:p-valuegroupscs}
\end{table}
\end{landscape}
}
\clearpage

\section{Centrality measures}{\label{app:centrality}}

\begin{table}[h!]
\resizebox{255pt}{!}{%
\begin{tabular}{l||llll}
     & \textit{Betweenness}    & \textit{Degree}    & \textit{Closeness}    & \textit{Local reaching}    \\ \hline\hline
LEMD & 0.0228 & 0.2736 & 0.1659 & 0.0227 \\
LEBL & 0.0176 & 0.2606 & 0.163  & 0.0175 \\
LEPA & 0.0113 & 0.1866 & 0.1457 & 0.0113 \\
LEMG & 0.0064 & 0.1589 & 0.1401 & 0.0064 \\
LEAL & 0.0043 & 0.124  & 0.1321 & 0.0043 \\
GCLP & 0.0037 & 0.0965 & 0.1267 & 0.0037 \\
LEIB & 0.0031 & 0.0961 & 0.1296 & 0.0031 \\
LEVC & 0.0019 & 0.0886 & 0.1286 & 0.0019 \\
GCTS & 0.0016 & 0.0844 & 0.1198 & 0.0016 \\
LEGE & 0.0016 & 0.0462 & 0.0958 & 0.0016 \\
LEZL & 0.0011 & 0.0677 & 0.1248 & 0.0011 \\
LEZG & 0.001  & 0.0203 & 0.0983 & 0.001  \\
GCRR & 0.0007 & 0.0557 & 0.1172 & 0.0007 \\
LEMH & 0.0007 & 0.0416 & 0.1192 & 0.0007 \\
GCXO & 0.0007 & 0.0249 & 0.1166 & 0.0007 \\
GCFV & 0.0006 & 0.0532 & 0.1164 & 0.0006 \\
LEBB & 0.0004 & 0.053  & 0.1256 & 0.0004 \\
LEST & 0.0003 & 0.0302 & 0.1196 & 0.0003 \\
LEVD & 0.0003 & 0.006  & 0.1005 & 0.0003 \\
LECU & 0.0002 & 0.0049 & 0.0731 & 0.0003 \\
LERS & 0.0002 & 0.0207 & 0.0876 & 0.0002 \\
LEJR & 0.0002 & 0.0195 & 0.1106 & 0.0002 \\
LEAM & 0.0002 & 0.0182 & 0.1108 & 0.0002 \\
LESO & 0.0002 & 0.0053 & 0.106  & 0.0002 \\
LEVX & 0.0002 & 0.0117 & 0.1084 & 0.0002 \\
LEGR & 0.0002 & 0.0123 & 0.1139 & 0.0002 \\
LEPP & 0.0002 & 0.0048 & 0.1012 & 0.0002 \\
LEAS & 0.0002 & 0.0185 & 0.1167 & 0.0002 \\
LEXJ & 0.0001 & 0.0154 & 0.1115 & 0.0001 \\
LEVT & 0.0001 & 0.0152 & 0.0962 & 0.0001 \\
LECO & 0.0001 & 0.0113 & 0.1096 & 0.0001 \\
LELL & 0.0001 & 0.0022 & 0.0303 & 0.0001 \\
LELN & 0.0001 & 0.0031 & 0.0842 & 0.0001 \\
LESA & 0.0001 & 0.0035 & 0.0563 & 0.0001 \\
LEAB & 0.0001 & 0.0011 & 0.0179 & 0.0001 \\
LEBZ & 0.0001 & 0.0032 & 0.0904 & 0.0001 \\
LESB & 0.0001 & 0.0014 & 0.022  & 0.0001 \\
LERJ & 0.0    & 0.0015 & 0.0709 & 0.0    \\
GCLA & 0.0    & 0.0101 & 0.1045 & 0.0    \\
GEML & 0.0    & 0.0039 & 0.1048 & 0.0    \\
LEBG & 0.0    & 0.0008 & 0.0374 & 0.0    \\
LEHC & 0.0    & 0.0004 & 0.0088 & 0.0    \\
LEBA & 0.0    & 0.0006 & 0.0191 & 0.0    \\
GCGM & 0.0    & 0.0021 & 0.08   & 0.0    \\
GCHI & 0.0    & 0.0033 & 0.0846 & 0.0   
\end{tabular}}
\caption{}
\label{table:centrality}
\end{table}

\clearpage
\nocite{*}
\printbibliography
\end{document}